\documentclass[11pt]{article}

\usepackage{amssymb}
\usepackage{epsfig}
\usepackage{hyperref}

\newtheorem{theorem}{Theorem}[section]
\newtheorem{proposition}[theorem]{Proposition}
\newtheorem{lemma}[theorem]{Lemma}

\newtheorem{example}[theorem]{Example}
\newtheorem{construction}[theorem]{Construction}

\newcommand{\hal}{\frac{1}{2}}
\newcommand{\Qed}{\rule{2.5mm}{3mm}}
\newcommand{\half}{$\frac{1}{2}$}
\newcommand{\hatr}{half-arc-transitive }
\newcommand{\ZZ}{\mathbb{Z}}
\newcommand{\rhogr}{\langle \rho \rangle}
\newcommand{\Aut}{\mathrm{Aut}}
\newcommand{\ot}{\leftarrow}
\newcommand{\iso}{\cong}

\newcommand{\X}{\mathcal{X}}
\newcommand{\Y}{\mathcal{Y}}
\newcommand{\Z}{\mathcal{Z}}
\newcommand{\C}{\mathcal{C}}
\newcommand{\HH}{\mathcal{H}}
\newcommand{\PP}{\mathcal{P}}
\newcommand{\la}{\langle}
\newcommand{\ra}{\rangle}
\newcommand{\setmin}{\backslash}
\newcommand{\mmod}[3]{#1 \equiv #2\, (\mathrm{mod}\ #3)}
\newenvironment{proof}{{\noindent \sc Proof:}}{\hfill $\Qed$}

\begin{document}

\begin{center}
{\bf\large ON QUARTIC HALF-ARC-TRANSITIVE METACIRCULANTS}
\end{center}

\medskip\noindent
\begin{center}
   ~Dragan Maru\v si\v c$\,^{a,b}$\addtocounter{footnote}{0}$^,$\footnotemark$^,$*
 ~and Primo\v z \v Sparl$^b$$^,$\addtocounter{footnote}{-1}\footnotemark

\bigskip
 
{\it {\small $^a$University of Primorska, Cankarjeva 5, 6000 Koper, Slovenia\\
$^b$IMFM, University of Ljubljana, Jadranska 19, 1111 Ljubljana, Slovenia}}
\end{center}

\addtocounter{footnote}{0}
\footnotetext{Supported in part by ``ARRS -- Agencija za znanost Republike Slovenije'', program no. P1-0285.

~* Corresponding author e-mail: ~dragan.marusic@guest.arnes.si}

\bigskip

%%%%%%%%%%%%%%%%%%%%%%%%%%%%%%%%%%%%%%%%%%
 %%%%%%%%%%%%%%%%%%%%%%%%%%%%%%%%%%%%%%%%%%
 %%%%%%%%%%      Abstract      %%%%%%%%%%%%
 %%%%%%%%%%%%%%%%%%%%%%%%%%%%%%%%%%%%%%%%%%
 %%%%%%%%%%%%%%%%%%%%%%%%%%%%%%%%%%%%%%%%%%

\hrule
\begin{abstract}
Following Alspach and Parsons, a {\em metacirculant graph} is a graph admitting a transitive group generated by 
two automorphisms $\rho$ and $\sigma$, where $\rho$ is $(m,n)$-semiregular for some integers $m \geq 1$, $n \geq 2$,
and where $\sigma$ normalizes $\rho$,
cyclically permuting the orbits of $\rho$ in such a way that $\sigma^m$ has at least one fixed vertex.
A {\em half-arc-transitive graph} is a vertex- and edge- but not arc-transitive graph. 
In this article quartic half-arc-transitive metacirculants are explored and their connection to the
so called tightly attached quartic half-arc-transitive graphs is explored. It is shown that
there are three essentially different possibilities for a quartic half-arc-transitive metacirculant which is not tightly attached
to exist. These graphs are extensively studied and some infinite families of such graphs are constructed.
\end{abstract}

\begin{quotation}\small
\end{quotation}
\hrule\medskip

%%%%%%%%%%%%%%%%%%%%%%%%%%%%%%%%%%%%%%%%%%%
%%%%%%%%%%%%%%%%%%%%%%%%%%%%%%%%%%%%%%%%%%%
%%%%%%%   Introduction    %%%%%%%%%%%%%%%%%
%%%%%%%%%%%%%%%%%%%%%%%%%%%%%%%%%%%%%%%%%%%
%%%%%%%%%%%%%%%%%%%%%%%%%%%%%%%%%%%%%%%%%%%

\section{Introductory and historic remarks}
\label{sec:intro}

\indent

Throughout this paper graphs are assumed to be finite 
and, unless stated otherwise, simple, connected and undirected
(but with an implicit orientation of the edges when appropriate).
For group-theoretic concepts not defined here we
refer the reader to \cite{BW,DM,W}, and for graph-theoretic terms not 
defined here we refer the reader to \cite{BM}.

Given a graph $X$ we let $V(X)$, $E(X)$, $A(X)$ and $\Aut X$ be
the vertex set, the edge set, the arc set and the automorphism
group of $X$, respectively. 
A graph $X$ is said to be {\em vertex-transitive}, 
{\em edge-transitive} and {\em arc-transitive} if its automorphism 
group $\Aut X$ acts transitively on $V(X)$, $E(X)$ and $A(X)$, respectively.
We say that $X$ is {\em half-arc-transitive} provided it is
vertex- and edge- but not arc-transitive.
More generally, by a {\em half arc-transitive} action of a 
subgroup $G \leq \Aut X$ on $X$ we mean a vertex- and edge- but not
arc-transitive action of $G$ on $X$.  
In this case we say that the graph $X$ is 
$(G,\hal)$-{\em arc-transitive}, and we say that the graph $X$ is
$(G,\hal,H)$-{\em arc-transitive} when it needs to be stressed that 
the vertex stabilizers $G_v$ (for $v \in V(X)$) are isomorphic to a 
particular subgroup $H \leq G$.
By a classical result of Tutte \cite[7.35, p.59]{T66},
a graph admitting a half-arc-transitive group action is necessarily of even valency.
A few years later Tutte's question as to the existence  of half-arc-transitive graphs
of a given even valency was answered by Bouwer \cite{B70} 
with a construction of a $2k$-valent half-arc-transitive graph for every $k\geq2$.
The smallest graph in Bouwer's family has 54 vertices and
valency 4. Doyle~\cite{PD76} and Holt~\cite{H81} independently found one with 27 vertices,
a graph that is now known to be the smallest half-arc-transitive graph~\cite{AMN94}.
 
Interest in the study of this class of graphs
reemerged about a decade later following a series of papers
dealing mainly with classification of certain restricted classes of such graphs 
as well as with various methods of constructions of new families of such graphs
%\cite{AMN94,AX94,MX97,TX94,TW89,Wa94,Xu92}; 
\cite{AMN94,AX94,TX94,TW89,Wa94,Xu92}; 
but see also the survey article 
\cite{DMsur} which covers the respective literature prior to 1998.
With some of the research emphasis shifting to questions concerning 
structural properties of half-arc-transitive graphs,
these graphs have remained an active topic of research to this day;
%see \cite{AX94,CM02,DAN04,DX98, FWZ07, LiS01,LLM04, DM98,DM05,MN98,MNgrt,MNjams, MPi99, PSxx,Wa94,Xu92,ZF06}.
see \cite{CM02,DAN04,DX98, FWZ07, LiS01,LLM04,MM99, DM98,DM05,MN98,MNgrt,MNjams, MPi99, MP99,MW00,PSxx, SW04,ZF06}.

In view of the fact that $4$ is the smallest admissible valency
for a half-arc-transitive graph, 
special attention has rightly been given to 
the study of quartic half-arc-transitive graphs.
One of the possible approaches
in the investigation of their properties 
concerns the so called  "attachment of alternating cycles" question.
Layed out in \cite{DM98}, the underlying theory
is made up of the following main ingredients.
For a quartic graph $X$ 
admitting a half-arc-transitive action of some subgroup
$G$ of $\Aut X$, let $D_G(X)$ be one of the two oriented
graphs associated in a natural way with the action of $G$ on $X$.
(In other words, $D_G(X)$ is an orbital graph
of $G$ relative to a non-self-paired orbital associated 
with a non-self-paired suborbit of length $2$
and $X$ is its underlying undirected graph.)
An even length cycle $C$ in $X$ is a
$G$-{\em alternating cycle} if every other vertex of $C$
is the tail and every other vertex of $C$ is the head
(in $D_G(X)$) of its two incident edges.
It was shown in \cite{DM98} that, first,
all $G$-alternating cycles of $X$ have the same length
-- half of this length is called the $G$-{\em radius} of $X$ --
and second, that any two adjacent $G$-alternating cycles
intersect in the same number of vertices, called
the $G$-{\em attachment number} of $X$. The intersection of two adjacent $G$-alternating cycles is called
a $G$-{\em attachment set}.
The attachment of alternating cycles concept has been addressed 
in a number of papers \cite{DM98,MP99,MW00,PSxx, SW04}
with a particular attention given
to the so called $G$-{\em tightly attached graphs},
that is, graphs where two adjacent $G$-alternating cycles have every other vertex in common.
In other words, their $G$-attachment number coincides with $G$-radius.
In all the above definitions the symbol $G$ is omitted when $G = \Aut X$.
Tightly attached graphs with odd radius have been completely classified 
in \cite{DM98}, whereas the classification of tightly attached graphs
with even radius, dealt with also in \cite{MP99,SW04}, has been
very recently completed in \cite{PSxx}.
At the other extreme, graphs with $G$-attachment number equal to $1$
and $2$, respectively, are called $G$-{\em loosely attached graphs}
and $G$-{\em antipodally attached graphs}.
As shown in \cite{MW00}, there exist infinite families 
of quartic half-arc-transitive graphs with arbitrarily prescribed attachment numbers. 
However, in view of the fact that every quartic half-arc-transitive
graph may be obtained as a cover of a loosely, antipodally or tightly attached graph \cite{MP99},
it is these three families of graphs that deserve special attention.

Now, as it turns out, all tightly attached quartic half-arc-transitive graphs
are metacirculant graphs \cite{DM98}. (For the definition of a metacirculant graph see
Section~\ref{sec:examples}.)  
The connection between the two classes of graphs 
goes so far as to suggest that even if quartic half-arc-transitive metacirculants which are 
not tightly attached do exist, constructing them will not be an easy task.
Exploring this connection is the main aim of this article.
Although short of a complete classification of quartic half-arc-transitive metacirculants,
we obtain a description of the three essentially different possibilities for a quartic half-arc-transitive metacirculant
which is not tightly attached to exist, together with constructions of infinite families of such graphs.
In doing so we give a natural decomposition of quartic half-arc-transitive metacirculants
into four classes depending on the structure of the quotient circulant graph relative to the semiregular automorphism $\rho$.
Loosely speaking, Class~I consists of those graphs whose quotient graph is a 'double-edged' cycle, Class~II
consists of graphs whose quotient is a cycle with a loop at each vertex, Class~III consists of graphs
whose quotient is a circulant of even order with antipodal vertices joined by a double edge, and Class~IV
consists of graphs whose quotient is a quartic circulant which is a simple graph (see Figure~\ref{fig:classes}).

The paper is organized as follows.
Section~\ref{sec:examples} contains some terminology together with four infinite families of quartic metacirculants, 
playing an essential role in the rest of the paper.
Section~\ref{sec:classes} gives the above mentioned decomposition.
Section~\ref{sec:ClassI} is devoted to Class~I graphs; in particular it is shown that this class coincides with the
class of tightly attached graphs (see Theorem~\ref{the:TA=2}).
Next, Section~\ref{sec:ClassII} deals with Class~II graphs. A characterization of the graphs of this class which are 
not tightly attached is given (see Theorem~\ref{the:IItheorem}) enabling us to construct an infinite family
of such graphs (see Construction~\ref{cons:infinite}). Moreover, 
a list of all quartic half-arc-transitive metacirculants of Class~II, of order at most $1000$,
that are not tightly attached is given.
Finally, in Section~\ref{sec:conclusions} a construction of an infinite family of 
loosely attached (and thus not tightly attached) half-arc-transitive metacirculants of Class~IV is given.

%%% ta Section dokoncno zlikan 26. 1. 2007

%%%%%%%%%%%%%%%%%%%%%%%%%%%%%%%%%%%%%%%%%%%
%%%%%%%%%%%%%%%%%%%%%%%%%%%%%%%%%%%%%%%%%%%
%%%%%%%  Definitions and examples    %%%%%%
%%%%%%%%%%%%%%%%%%%%%%%%%%%%%%%%%%%%%%%%%%%
%%%%%%%%%%%%%%%%%%%%%%%%%%%%%%%%%%%%%%%%%%%

\section{Definitions and examples}
\label{sec:examples}

\indent

We start by some notational conventions used throughout this paper.
Let $X$ be a graph. The fact that $u$ and $v$ are adjacent vertices of $X$ will be denoted by $u \sim v$; 
the corresponding edge will be denoted by $[u,v]$,
in short by $uv$.
In an oriented graph the fact that the edge $uv$ is oriented from $u$ to $v$
will be denoted by $u \to v$ (as well as by $v \ot u$). In this case 
the vertex $u$ is referred to as the {\em tail} and $v$ is referred to as the {\em head} of the edge $uv$. 
Let $U$ and $W$ be disjoint subsets of $V(X)$.
The subgraph of $X$ induced by $U$ will be denoted by $X[U]$;
in short, by $[U]$, when the graph $X$ is clear from the context.
Similarly, we let $X[U,W]$ (in short $[U,W]$) denote the bipartite subgraph 
of $X$ induced by the edges having one endvertex in $U$ and the other endvertex in $W$. 
Furthermore, if $\rho$ is an automorphism of $X$, we denote the corresponding {\em quotient
(multi)graph relative to $\rho$}, whose
vertex set is the set of orbits of $\rho$ with two orbits adjacent whenever there is an edge in $X$ 
joining vertices from these two orbits, by $X_\rho$. 

For the sake of completeness we include the definition of a Cayley graph. 
Given a group $G$ and an inverse closed subset $S \subseteq G \setmin \{1\}$ the
{\em Cayley graph} $\mathrm{Cay}(G\, |\, S)$ is the graph with vertex set $G$ and edges of 
the form $[g, gs]$, where $g \in G$, $s \in S$.

Let $m \geq 1$ and $n \geq 2$ be integers.
An automorphism of a graph is called $(m,n)$-{\em semiregular}
if it has $m$ orbits of length $n$ and no other orbit. 
We say that a graph $X$ is an $(m,n)$-{\em metacirculant graph} (in short an 
$(m,n)$-{\em metacirculant}) 
if there exists an $(m,n)$-semiregular automorphism $\rho$ of $X$, 
together with an additional automorphism
$\sigma$ of $X$ normalizing $\rho$, that is,
\begin{equation} \label{eq:metagroup}
	\sigma^{-1}\rho\sigma = \rho^r\quad \mathrm{for}\ \mathrm{some}\quad r \in
	\mathbb{Z}_n^*,
\end{equation}
and cyclically permuting the orbits of $\rho $ in such a way that 
$\sigma^m$ fixes a vertex of $X$. (Hereafter $\ZZ_n$ denotes the ring of residue classes
modulo $n$ as well as the additive cyclic group of order $n$, depending on the context.) Note that this implies
that $\sigma^m$ fixes a vertex in every orbit of $\rho$. 
To stress the role of these two automorphisms in the definition of 
the metacirculant $X$ we shall say that $X$ 
is an $(m,n)$-{\em metacirculant relative to the ordered pair} $(\rho,\sigma)$.
Obviously, a graph is an $(m,n)$-metacirculant
relative to more than just one ordered pair of automorphisms
except for the trivial case when $m=1$ and $n=2$, which  corresponds
to $X \cong K_2$.
For example, the automorphism $\sigma$ may be replaced by 
$\sigma\rho$.
A graph $X$ is a {\em metacirculant} if it is an $(m,n)$-metacirculant for
some $m$ and $n$.
This definition is equivalent with the original definition of a metacirculant
by Alspach and Parsons (see \cite{AP82}).
For the purposes of this paper we extend this definition somewhat.
We say that a graph $X$ is a {\em weak} $(m,n)$-{\em metacirculant} 
(more precisely a {\em weak} $(m,n)$-{\em metacirculant relative to the 
ordered pair $(\rho,\sigma)$})
if it has all the properties of an $(m,n)$-metacirculant except that we do not
require that $\sigma^m$ fixes a vertex of $X$.
We say that $X$ is a {\em weak metacirculant} if it is
a weak $(m,n)$-metacirculant for some positive integers $m$ and $n$.

Note that there exist integers $m$, $n$ and
weak $(m,n)$-metacirculants which are not $(m,n)$-metacirculants.
For example, the graph $\Y(10,100;11,90)$ 
(see Example~\ref{ex:Y} below)
is a weak $(10,100)$-metacirculant but it can be seen that it
is not a $(10,100)$-metacirculant.
However, this graph is also a $(40,25)$-metacirculant. 
The question remains if the class of weak metacirculants is indeed 
larger than that of metacirculants.
Nevertheless, at least for the purposes of this paper it proves natural
to work in the context of weak metacirculants.

%%%%%%%%%%%%%%%%%%%%%%%%%%%%%%%
%%%%%%%% EXAMPLES %%%%%%%%%%%%%
%%%%%%%%%%%%%%%%%%%%%%%%%%%%%%%
Below we give a few infinite families of weak metacirculants that will play
a crucial role in the investigation of quartic half-arc-transitive metacirculants, 
the main theme of this article.

%%%%%%  ex:Xodd  %%%%%%
\begin{example}
\label{ex:Xodd}\rm
For each $m \geq 3$, for each odd $n \geq 3$ and for each  $r \in \ZZ_n^*$, where $r^m = \pm 1$, 
let $\X_o(m,n;r)$ be the graph with vertex 
set $V = \{u_i^j\ |\ i \in \ZZ_m,\ j \in \ZZ_n\}$ 
and edges defined by the following adjacencies:
$$ u_i^j \sim u_{i+1}^{j \pm r^i}\quad ;\quad i \in \ZZ_m,\ j \in \ZZ_n.$$
(Note that the subscript $o$ in the symbol $\X_o(m,n;r)$ is meant to 
indicate that $n$ is an odd integer.)
The permutations $\rho$ and $\sigma$, defined 
by the rules
$$ \begin{array}{lll} u_i^j\rho = u_i^{j+1} &  ; & i \in \ZZ_m,\  j \in \ZZ_n \\ \\
     u_i^j\sigma = u_{i+1}^{rj} & ; & i \in \ZZ_m,\ j \in \ZZ_n, \end{array}$$ 
are automorphisms of $\X_o(m,n;r)$. Note that $\rho$ is $(m,n)$-semiregular and that
$\sigma^{-1}\rho\sigma = \rho^r$. Moreover, $\sigma$ cyclically permutes the orbits of 
$\rho$ and $\sigma^m$ fixes $u_i^0$ for every $i \in \ZZ_m$. Hence $\X_o(m,n;r)$
is an $(m,n)$-metacirculant. 
We note that graphs $\X_o(m,n;r)$ correspond to the graphs $X(r;m,n)$ introduced in \cite{DM98}.
We also note that the Holt graph, 
the smallest half-arc-transitive graph (see \cite{AMN94,PD76,H81}), 
is isomorphic to $\X_o(3,9;2)$. 
\end{example}

%%%%%%  ex:Xeven  %%%%%%
\begin{example}
\label{ex:Xeven}\rm
For each $m \geq 4$ even, $n \geq 4$ even, $r \in \ZZ_n^*$, where $r^m = 1$, and $t \in \ZZ_n$,
where $t(r-1) = 0$,
let $\X_e(m,n;r,t)$ be the graph with vertex 
set $V = \{u_i^j\ |\ i \in \ZZ_m,\ j \in \ZZ_n\}$ 
and edges defined by the following adjacencies:
$$ u_i^j \sim \left\{\begin{array}{lll}
	u_{i+1}^j,\ u_{i+1}^{j + r^i} & ; & i \in \ZZ_m \setmin \{m-1\},\ j \in \ZZ_n \\ \\
	u_{0}^{j+t},\ u_0^{j+r^{m-1}+t} & ; & i = m-1,\ j \in \ZZ_n .\end{array}\right. $$
(In analogy with Example~\ref{ex:Xodd} the subscript $e$ in the symbol $\X_e(m,n;r,t)$ is meant to 
indicate that $n$ is an even integer.)
The permutations $\rho$ and $\sigma$, defined 
by the rules
$$ 	u_i^j\rho = u_i^{j+1} \quad  ; \quad i \in \ZZ_m,\  j \in \ZZ_n $$
$$	u_i^j\sigma = \left\{\begin{array}{lll}
	u_{i+1}^{rj} & ; & i \in \ZZ_m \setmin \{m-1\},\ j \in \ZZ_n\\ \\
	u_{0}^{rj+t} & ; & i = m-1,\ j \in \ZZ_n ,\end{array}\right.$$ 
are automorphisms of $\X_e(m,n;r,t)$. Note that $\rho$ is $(m,n)$-semiregular, that
$\sigma^{-1}\rho\sigma = \rho^r$ and that $\sigma$ cyclically permutes the orbits of 
$\rho$. Hence $\X_e(m,n;r,t)$ is a weak $(m,n)$-metacirculant. 
\end{example}

As noted in Section~\ref{sec:intro}, a complete classification of quartic tightly attached 
half-arc-transitive graphs is given in
\cite{DM98} for odd radius and in \cite{PSxx} for even radius. 
It follows by this classification that a quartic
tightly attached half-arc-transitive graph is isomorphic either to some $\X_o(m,n;r)$ 
or to some $\X_e(m,n;r,t)$, depending on the radius parity.

%%%%%%  ex:Y  %%%%%%
\begin{example}
\label{ex:Y}\rm
For each $m \geq 3$, $n \geq 3$, $r \in \ZZ_n^*$, where $r^m = 1$, and $t \in \ZZ_n$ satisfying $t(r-1)=0$,
let $\Y(m,n;r,t)$ be the graph with vertex 
set $V = \{u_i^j\ |\ i \in \ZZ_m,\ j \in \ZZ_n\}$ 
and edges defined by the following adjacencies:
$$ u_i^j \sim \left\{\begin{array}{lll}
	u_i^{j + r^i},\ u_{i + 1}^j & ; & i \in \ZZ_m \setmin \{m-1\},\ j \in \ZZ_n \\ \\
	u_{m-1}^{j + r^{m-1}},\ u_0^{j+t} & ; & i = m-1,\ j \in \ZZ_n .\end{array}\right. $$
The permutations $\rho$ and $\sigma$, defined 
by the rules
$$ u_i^j\rho = u_i^{j+1} \quad ; \quad  i \in \ZZ_m,\  j \in \ZZ_n $$
$$ u_i^j\sigma = \left\{\begin{array}{lcl}
		u_{i+1}^{rj} & ; & i \in \ZZ_m \setmin \{m-1\},\ j \in \ZZ_n\\ \\
		u_0^{rj+t} & ; & i = m-1,\ j \in \ZZ_n ,\end{array}\right.$$
are automorphisms of $\Y(m,n;r,t)$. Observe that $\rho$ is $(m,n)$-semiregular and that
$\sigma^{-1}\rho\sigma = \rho^r$. Moreover, $\sigma$ cyclically permutes the orbits of $\rho$, and so 
$\Y(m,n;r,t)$ is a weak $(m,n)$-metacirculant.
We note that the Holt graph, see Example~\ref{ex:Xodd}, is also isomorphic to $\Y(3,9;7,3)$.
\end{example}

%%%%%%  ex:Z  %%%%%%
\begin{example}
\label{ex:Z}\rm
For each $m \geq 5$, $n \geq 3$, $k \in \ZZ_m \setmin \{0,1,-1\}$ and $r \in \ZZ_n^*$, where $r^m = 1$,
let $\Z(m,n;k,r)$ be the graph with vertex 
set $V = \{u_i^j\ |\ i \in \ZZ_m,\ j \in \ZZ_n\}$ 
and edges defined by the following adjacencies:
$$ u_i^j \sim u_{i+1}^j,\ u_{i+k}^{j+r^i}\quad ; \quad i \in \ZZ_{m},\ j\in \ZZ_n .$$
The permutations $\rho$ and $\sigma$, defined 
by the rules
$$ u_i^j\rho = u_i^{j+1} \quad ; \quad  i \in \ZZ_m,\  j \in \ZZ_n $$
$$ u_i^j\sigma = u_{i+1}^{rj} \quad ; \quad i \in \ZZ_m,\  j \in \ZZ_n ,$$
are automorphisms of $\Z(m,n;k,r)$. Observe that $\rho$ is $(m,n)$-semiregular and that 
$\sigma^{-1}\rho\sigma = \rho^r$. Moreover, $\sigma$ cyclically permutes the orbits of 
$\rho$ and $\sigma^m = 1$. Hence $\Z(m,n;k,r)$ is an $(m,n)$-metacirculant.
\end{example}

%%%%%%%%%%%%%%%%%%%%%%%%%%%%%%%%%%%%%%%%%%%
%%%%%%%%%%%%%%%%%%%%%%%%%%%%%%%%%%%%%%%%%%%
%%   The four classes %%
%%%%%%%%%%%%%%%%%%%%%%%%%%%%%%%%%%%%%%%%%%%
%%%%%%%%%%%%%%%%%%%%%%%%%%%%%%%%%%%%%%%%%%%

\section{The four classes}
\label{sec:classes}

\indent

In this section we start our investigation
of half-arc-transitivity of quartic weak metacirculants.
First, we state a result from \cite{DM98} which will be used throughout the rest of the paper.

\begin{proposition}[{\cite[Proposition~2.1.]{DM98}}] \label{pro:flip}
Let $X$ be a half-arc-transitive graph. Then no automorphism of $X$ can
interchange a pair of adjacent vertices in $X$.
\end{proposition}

Throughout this section we let $X$ denote
a connected quartic half-arc-transitive weak $(m,n)$-metacirculant 
relative to an ordered pair $(\rho,\sigma)$. Furthermore,
we let $X_i$, $i \in \ZZ_m$, denote the orbits of $\rho$ where 
$X_{i+1} = X_i\sigma$ for each $i \in \ZZ_m$.
Clearly, the degrees of subgraphs $[X_i]$ are all equal. 
We shall denote this number by
$d_{inn}(X)$ and call it the {\em inner degree} of $X$. 
Note that $d_{inn}(X)$ must be even, for otherwise $n$ is even and a vertex $u$
of $X$ is necessarily adjacent to $u\rho^{\frac{n}{2}}$. 
But then $\rho^{\frac{n}{2}}$ interchanges two adjacent vertices,
which contradicts Proposition~\ref{pro:flip}.
Furthermore, $d_{inn}(X)$ cannot be $4$, for otherwise the connectedness of $X$ implies 
that $m = 1$ and thus
$X$ is a circulant. But no \hatr Cayley graph of an abelian group exists. 
Namely, if $X = \mathrm{Cay}(G\, |\, S)$ choose $s \in S$, let $\varphi : X \to X$ map $x$ to $x^{-1}$ and let
$\rho_{s} : X \to X$ map $x$ to $sx$. Then $\varphi$ and $\rho_s$ are automorphisms of $X$ and 
$\varphi\rho_s$ interchanges adjacent vertices $1$ and $s$.
Therefore
\begin{equation}
\label{eq:d_inn}
	d_{inn}(X) \in \{0,2\}.
\end{equation}

We now show that the number of orbits of $\rho$ is at least $3$.

%%%%%% pro:threeorbits  %%%%%%
\begin{proposition}
\label{pro:threeorbits}
Let $X$ be a connected quartic \hatr weak $(m,n)$-metacirculant relative to
an ordered pair $(\rho, \sigma)$.
Then $m \geq 3$.
\end{proposition}

\begin{proof}
By the above remarks we have $d_{inn}(X) \in \{0,2\}$ and $m \geq 2$. 
Assume then that $m =2$ and let $U$ and $W$ be the orbits of $\rho$.
We show that there exists an automorphism of $X$ fixing  
$U$ and $W$ setwise and interchanging two adjacent vertices of $U$ which
contradicts Proposition~\ref{pro:flip}. By
\cite[Proposition~2.2.]{DM98}, which states that a graph cannot be half-arc-transitive if it has 
a $(2,n)$-semiregular automorphism whose two orbits give rise to a bipartition of the graph in question, we 
must have $d_{inn}(X) = 2$. 
Fix a vertex $u \in U$ and set $u^i = u\rho^i$, where $i \in \ZZ_n$. 
There exists some nonzero
$s \in \mathbb{Z}_n$ such that $u^i \sim u^{i\pm s}$ for all $i \in \ZZ_n$.
Next, choose a vertex $w \in W$ such that $u^0 \sim w$ and set
$w^i = w\rho^i$, where $i \in \ZZ_n$. Letting $r \in \ZZ_n^*$ be as in
equation (\ref{eq:metagroup}) we have $u^0\sigma \sim u^{\pm s}\sigma =
u^0\rho^{\pm s}\sigma = u^0\sigma \rho^{\pm rs}$, and so
$w^i \sim w^{i\pm rs}$ for all $i \in \ZZ_n$.
There exists some nonzero $t \in \mathbb{Z}_n$ such that
$u^0 \sim w^t$. Therefore, we have $u^i \sim w^i,w^{i+t}$ for all $i \in \ZZ_n$.
It is easy to see that the permutation $\varphi$ of $V(X)$ defined by the rule $u^i\varphi = u^{-i}$ and
$w^i\varphi = w^{t-i}$, where $i \in \ZZ_n$, is an automorphism
of $X$. But then $\varphi\rho^s$ interchanges adjacent vertices
$u^0$ and $u^{s}$, completing the proof of Proposition~\ref{pro:threeorbits}.
\end{proof}
\bigskip

\begin{figure}
\begin{center}
\includegraphics[scale=0.35]{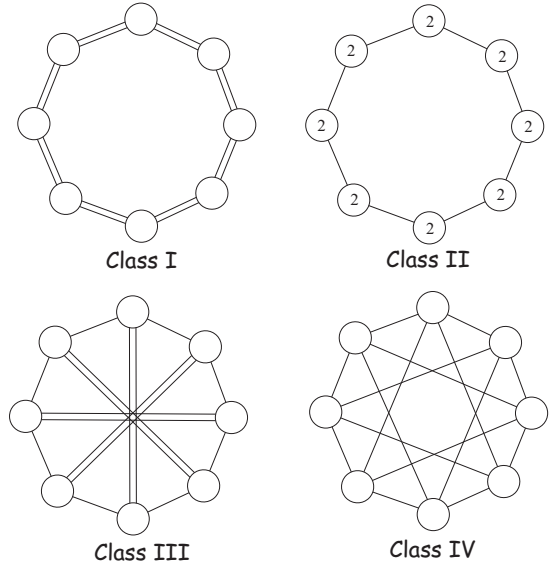}
\caption{Every quartic half-arc-transitive metacirculant falls
into one (or more) of the four classes.}
\label{fig:classes}
\end{center}
\end{figure}

We now use (\ref{eq:d_inn}) and Proposition~\ref{pro:threeorbits} 
to show that each connected quartic half-arc-transitive weak metacirculant
belongs to at least one of the following four  classes reflecting
four essentially different ways in which a quartic graph may be
a half-arc-transitive weak metacirculant (see Figure~\ref{fig:classes}). These four classes are
described below. (Recall that the orbits of $\rho$ are denoted by $X_i$.)

\begin{itemize}
   \item {\bf Class I.} The graph $X$ belongs to {\em Class~I} if 
   $d_{inn}(X) = 0$ and each orbit $X_i$ is connected (with a double edge) 
   to two other orbits. In view of connectedness of $X$, we have that
   $X_\rho$ is a "double-edge" cycle.
   \item {\bf Class II.} The graph $X$ belongs to {\em Class~II} if 
   $d_{inn}(X) = 2$ and each orbit $X_i$ is connected (with a single edge) to 
   two other orbits. In view of connectedness of $X$, we have that
   $X_\rho$ is a cycle (with a loop at each vertex).
   \item {\bf Class III.} The graph $X$ belongs to {\em Class~III} if
   $d_{inn}(X) = 0$ and each orbit $X_i$
   is connected to three other orbits, to one with a double edge and
   to two with a single edge. Clearly, $m$ must be even in this case and 
   an orbit $X_i$ is connected to the orbit $X_{i + \frac{m}{2}}$ with a double edge.
   In short, $X_\rho$ is a connected circulant with double edges connecting antipodal vertices.
   \item {\bf Class IV.} The graph $X$ belongs to {\em Class~IV} if 
   $d_{inn}(X) = 0$ and each orbit $X_i$ is connected (with a single edge) to
   four other orbits. In short, $X_\rho$ is a connected circulant of valency $4$ and is a simple graph.
\end{itemize}

We remark that these four classes of metacirculants are not disjoint. For instance, it may be seen that
the Holt graph $\X_o(3,9;2) \cong \Y(3,9;7,3)$ belongs to Classes I and II but not to Classes III and IV.
Its canonical double cover, the smallest example in the Bouwer's construction, belongs to Classes I, II and III but not to Class IV.
On the other hand, the graph $\Z(20,5;9,2)$ belongs solely to Class IV. 

In the next two sections Classes I and II are analyzed in detail.  In the last section 
future research directions regarding
interconnectedness of Classes I, II, III and IV, are layed out.

%
%
% ta section  popravljen, December 17, 2006
%
%

%%%%%%%%%%%%%%%%%%%%%%%%%%%%%%%%%%%%%%%%%%%
%%%%%%%%%%%%%%%%%%%%%%%%%%%%%%%%%%%%%%%%%%%
%%%%%%%%%%%%%%%  Class I   %%%%%%%%%%%%%%%%
%%%%%%%%%%%%%%%%%%%%%%%%%%%%%%%%%%%%%%%%%%%
%%%%%%%%%%%%%%%%%%%%%%%%%%%%%%%%%%%%%%%%%%%

\section{Graphs of Class~I}
\label{sec:ClassI}

\indent

The aim of this section is to prove the following theorem.

\begin{theorem}
\label{the:TA=2}
Connected quartic half-arc-transitive weak metacirculants of Class~I coincide with 
connected quartic tightly attached half-arc-transitive graphs.
\end{theorem}

Throughout this section we let $X$ denote a connected 
quartic \hatr
weak metacirculant of Class~I and we let $m,n, \rho$ and $\sigma$ be such that 
$X$ is a weak $(m,n)$-metacirculant relative to
the ordered pair $(\rho,\sigma)$. Fix a vertex $u \in V(X)$ and let
$u_i^0 = u\sigma^i$ for all $i \in \{0,1,\ldots , m-1\}$. Then let
$u_i^j = u_i^0\rho^j$ for all $i \in \ZZ_m$, $j \in \ZZ_n$.
With this notation the orbits of $\rho$ are 
precisely the sets $X_i = \{u_i^j\ |\ j \in \ZZ_n \}$,
$i \in \ZZ_m$.
Since $X$ is connected we can assume that $X_0 \sim X_1$ in the
quotient graph $X_\rho$.
Moreover, as $X$ is a weak $(m,n)$-metacirculant relative to the pair 
$(\rho,\sigma\rho^j)$ for any $j \in \ZZ_n$, we can in fact assume
that $u_0^0$ is adjacent to $u_1^0$.
There exists some nonzero $a \in \ZZ_n$ such that
$u_0^0$ is adjacent also to $u_1^a$. 
Therefore $N(u_0^0) \cap X_1 = \{u_1^0, u_1^{a}\}$. Let $r \in \ZZ_n^*$ be as in equation
(\ref{eq:metagroup}). Then $u_1^a\sigma^i = u_1^0\rho^a\sigma^i = u_1^0\sigma^i\rho^{r^i a}$
holds for all $i \in \ZZ_m$ , and so
%%%%%%  eq:neighI  %%%%%%
\begin{equation}
\label{eq:neighI}
	N(u_i^j) \cap X_{i+1} = \{u_{i+1}^j, u_{i+1}^{j+r^i a}\} \quad \mathrm{for}\ \mathrm{all}\ 
		i \in Z_m \setmin \{m-1\},\ j \in \ZZ_n .
\end{equation}
Out of the two orientations of the edges of $X$ induced by the half-arc-transitive
action of $\Aut X$ we choose the one where $u_0^0 \to u_1^0$. Denote the corresponding 
oriented graph by $D_X$. There are two possibilities depending on whether 
$u_0^0$ is the tail or the head of the edge $u_0^0 u_1^a$ in $D_X$. 
In Lemma~\ref{le:class1++} below we show that in the former case $X$ is tightly attached, and 
in Lemma~\ref{le:class1+-} we show that the latter actually never occurs.

%%%%%%  le:class1++  %%%%%%
\begin{lemma} 
\label{le:class1++}
With the notation introduced in the previous paragraph, if $u_0^0$ is the tail of 
the edge $u_0^0 u_1^a$ in $D_X$ then $X$ is tightly attached.
\end{lemma}

\begin{proof}
Clearly in this case all the edges in $D_X[X_i, X_{i+1}]$ 
are oriented from $X_i$ to $X_{i+1}$. Therefore
(\ref{eq:neighI}) implies that $u_0^j \to u_1^j, u_1^{j + a}$ and that 
$u_1^j \to u_2^j, u_2^{j+ra}$ for all $j \in \ZZ_n$. 
Moreover, any alternating cycle of $X$ is a subgraph of $X[X_i,X_{i+1}]$ for
some $i \in \ZZ_m$.
We now inspect the two alternating cycles containing $u_1^0$. 
Denote the one containing vertices from $X_0$ and $X_1$ with $C_1$ and 
the one containing vertices from $X_1$ and $X_2$ with $C_2$.
In view of Proposition~\ref{pro:threeorbits} we have $m \geq 3$, and so
$C_1 \cap C_2 \subseteq X_1$. We have $C_1 \cap X_1 = \{u_1^j\, |\, j \in \langle a \rangle\}$. 
(Here $\la a \ra$ denotes the additive subgroup of $\ZZ_n$ generated by $a$.)
Moreover, $C_2 \cap X_1 = \{u_1^j\, |\, j \in \langle ra \rangle\}$.
But $r \in \ZZ_n^*$ and so  $\langle a \rangle = \langle ra \rangle$.
Hence $C_1 \cap X_1 = C_2 \cap X_2$, which completes the proof.
\end{proof}

%%%%%% le:class1+-  %%%%%%
\begin{lemma}
\label{le:class1+-}
There exists no connected quartic \hatr weak meta\-circulant of
Class~I such that, with the notation from the paragraph preceding the statement
of Lemma~\ref{le:class1++},  the vertex  $u_0^0$ is the
head of the edge $u_0^0 u_1^a$ in $D_X$.
\end{lemma}

\begin{proof}
Suppose that there does exist such a graph and denote it by $X$.
Our approach is as follows. We first show that the
stabilizer of a vertex in $X$ cannot be $\ZZ_2$. We then show that this forces
$m$ to be odd and $\mmod{n}{2}{4}$, which enables us to investigate the $\Aut X$-orbit
of the so called generic $8$-cycles of $X$ in a greater detail. In particular we find that $(r-1)^2 = 0$.
Then, investigating the $\Aut X$-orbit of a particular nongeneric
$8$-cycle, we finally arrive at a contradiction, thus showing that $X$ cannot exist.

Observe that since $\sigma^m$ fixes the orbits $X_i$ 
setwise, there exists some $t \in \ZZ_n$, such that $\sigma^m\rho^{-t}$ fixes $u_0^0$. 
Since the sets $X_i$ are blocks of imprimitivity for $\langle \rho, \sigma \rangle$, 
the particular orientation of the edges in $D_X$ implies that $\sigma^m\rho^{-t}$ fixes the neighbors of $u_0^0$
pointwise. Continuing this way we have, in view of connectedness of  $X$, that  $\sigma^m = \rho^t$. 
Note that equation (\ref{eq:metagroup}) implies that $\sigma^{-m}\rho\sigma^m = \rho^{r^m}$, and  so
%%%%%%  eq:leI+-r^m  %%%%%%
\begin{equation}
\label{eq:leI+-r^m}
	r^m = 1. 
\end{equation}
Moreover, (\ref{eq:metagroup}) also implies that $u_i^j\sigma = u_i^0\rho^j\sigma = u_i^0\sigma\rho^{rj}$, and  so
%%%%%%  eq:leI+-sigma  %%%%%%
\begin{equation}
\label{eq:leI+-sigma}
u_i^j\sigma = \left\{\begin{array}{lll}
	u_{i+1}^{rj} 	& ; 	& i \in \ZZ_m \setmin \{m-1\},\ j \in \ZZ_n \\ \\
	u_0^{rj+t} 	& ;	& i = m-1,\ j\in \ZZ_n.\end{array}\right.  
\end{equation}
Combining together (\ref{eq:neighI}) and  (\ref{eq:leI+-sigma}), we have that
for any $i \in \ZZ_m$ and $j \in \ZZ_n$, the two edges 
connecting $u_i^j$ to vertices from $X_{i+1}$ in $D_X$ are given by
%%%%%%  eq:leI+-adj  %%%%%%
\begin{equation}
\label{eq:leI+-adj}
	u_i^j \to \left\{\begin{array}{lll}
		u_{i+1}^{j} 	& ; 	& i \neq m-1 \\ \\
		u_0^{j+t} 	& ; 	& i = m-1\end{array}\right.
  	\mathrm{and}\quad 
	u_i^j \ot \left\{\begin{array}{lll}
		u_{i+1}^{j+r^i a} 	& ; 	& i \neq m-1 \\ \\
		u_0^{j+r^{m-1} a+t} 	& ;	& i = m-1. \end{array}\right.
\end{equation}
Since $\sigma$ maps the edge $u_{m-1}^0 u_0^t$ to the edge $u_0^t u_1^{rt}$, we also have 
%%%%%%  eq:leI+-rt  %%%%%%
\begin{equation}
\label{eq:leI+-rt}
	t(r-1) = 0.
\end{equation}

\noindent
{\sc Claim 1:} We lose no generality in assuming that $a = 1$.
\smallskip

Observe that since $X$ is connected,  (\ref{eq:leI+-adj}) implies that
$\langle a, t \rangle = \ZZ_n$. 
Let $d = |a|$ denote the order of $a$ in the additive  group $\ZZ_n$.
Clearly Claim~1 holds if $d = n$, so assume that $d < n$ and set $k = \frac{n}{d}$.
Let $\rho' = \rho^a$. Then $\rho'$ is an $(mk,d)$-semiregular automorphism of $X$ and
$\sigma^{-1}\rho'\sigma = \rho^{ra} = \rho'^r$. We now show that $\sigma$ 
cyclically permutes the orbits of $\rho'$. The orbit of $\rho'$ containing $u_0^0$ 
is $X_0' = \{u_0^{ja}\ |\ j \in \ZZ_n\}$. Moreover, $X_i' = X_0'\sigma^i = 
\{u_i^{jr^i a}\ |\ j \in \ZZ_n\}$ for $i = 0,1, \ldots , m-1$ and
$X_m' = X_0'\sigma^m = \{u_0^{ja + t}\ |\ j \in \ZZ_n\}$.
Since $\langle a, t \rangle = \ZZ_n$ and $\langle a \rangle \neq \ZZ_n$, the set $X_m'$ 
(which is clearly an orbit of $\rho'$) cannot be equal to $X_0'$. Continuing this way we see
that $\sigma$ cyclically permutes the $mk$ orbits of $\rho'$. Thus $X$ is a weak $(mk,d)$-metacirculant 
of Class~I relative to the ordered pair $(\rho', \sigma)$. It is now clear that in the notation
of vertices of $X$ relative to the ordered pair $(\rho', \sigma)$ the corresponding parameter $a'$ is equal to $1$.
From now on we can therefore assume that $a = 1$.
\medskip

\noindent
{\sc Claim 2:} Let $v \in V(X)$. Then $|(\Aut X)_v| > 2$.
\smallskip

Suppose on the contrary that $(\Aut X)_v \cong \ZZ_2$ for some (and hence any) $v \in V(X)$.
Let $\tau$ be the unique nontrivial automorphism of $(\Aut X)_{u_0^0}$. Then  $\tau$ interchanges $u_1^0$ and
$u_{m-1}^{-r^{m-1}-t}$ and also interchanges $u_1^1$ and $u_{m-1}^{-t}$. 
We now determine the action of $\tau$ on the vertices of $X$ recursively as follows.
Since $\tau \notin \langle \rho, \sigma \rangle$ and since
$u_1^1\tau = u_{m-1}^{-t} = u_1^1\sigma^{m-2}\rho^{-r^{m-2}-t}$ we have 
(recall that  $(\Aut X)_{u_1^1} \cong \ZZ_2$) 
that $u_0^1\tau \neq u_0^1\sigma^{m-2}\rho^{-r^{m-2}-t} = u_{m-2}^{-t}$. It follows that
$\tau$ interchanges $u_0^1$ and $u_0^{r^{m-1}}$. Now $u_0^1 \ot u_1^2$ and a similar argument 
shows that $\tau$ interchanges $u_1^2$ and $u_{m-1}^{r^{m-1}-t}$. 
Continuing this way we find that $\tau$ maps according to the rule:

%%%%%%  eq:leI+-tau
\begin{equation}
\label{eq:leI+-tau}
u_i^j\tau = \left\{\begin{array}{lll}
		u_{m-i}^{jr^{m-1} - r^{m-1}- r^{m-2} - \cdots -r^{m-i} -t} & 
			; 	& i \in \ZZ_m \setmin \{0\},\ j \in \ZZ_n\\ \\
		u_0^{jr^{m-1}}	& ; 	& i = 0,\ j \in \ZZ_n .\end{array}\right.
\end{equation}
We leave the details to the reader.
Recall now that, by assumption,  $\tau$ interchanges  $u_{m-1}^{-t}$ and
$u_1^1$. On the other hand,  (\ref{eq:leI+-tau}) implies that
$u_{m-1}^{-t}\tau = u_1^{-tr^{m-1} - r^{m-1} - r^{m-2} - \cdots - r - t}$. 
Therefore, equation (\ref{eq:leI+-rt}) implies that
%%%%%%  eq:leI+-cond  %%%%%%
\begin{equation}
\label{eq:leI+-cond}
1+r+r^2+\cdots + r^{m-1} + 2t = 0.
\end{equation}
We now define a mapping $\psi$ on $V(X)$ by the rule
%%%%%%  eq:leI+-psi
\begin{equation}
\label{eq:leI+-psi}
u_i^j\psi = \left\{\begin{array}{lll}
		u_i^{-j+1+r+r^2+\cdots + r^{i-1}} & ; 	& i \in \ZZ_m \setmin \{ 0 \},\ j \in \ZZ_n \\ \\
		u_0^{-j}	& ; 	& i = 0,\ j \in \ZZ_n .\end{array}\right. 
\end{equation}
Clearly $\psi$ is a bijection. It is easy to check that $\psi$ maps every edge
of $[X_i,X_{i+1}]$, where $i \in \ZZ_m \setmin \{m-1\}$, to an edge of $X$.
As for the edges of $[X_{m-1},X_0]$, note that by (\ref{eq:leI+-adj}) we have 
that $N(u_{m-1}^j) \cap X_0 = \{u_0^{j+t}, u_0^{j+r^{m-1}+t}\}$.
Observe that $\psi$ maps the latter two vertices to 
$u_0^{-j-t}$ and $u_0^{-j-r^{m-1}-t}$, respectively. Moreover, in view of
equation (\ref{eq:leI+-cond}) we have that
$u_{m-1}^j \varphi = u_{m-1}^{-j+1+r+r^2+\cdots + r^{m-2}} = u_{m-1}^{-j-r^{m-1}-2t}$, 
and so $\psi$ is an automorphism of $X$. 
Since $X$ is half-arc-transitive, there exists some $\varphi \in \Aut X$ mapping the
edge $u_1^1 u_0^0$ of $D_X$ to the edge $u_0^0 u_1^0$. But then $\psi \varphi$ interchanges
adjacent vertices $u_0^0$ and $u_1^0$,  contradicting Proposition~\ref{pro:flip}.
Therefore, $|(\Aut X)_v| > 2$, as claimed.
\medskip

Let now $A$ be the attachment set of $X$ containing $u_0^0$. 
In view of \cite[Lemma~3.5.]{MP99}, which states that
in a finite connected quartic \hatr graph with attachment sets containing at least three vertices, 
the vertex stabilizers are isomorphic to $\ZZ_2$, it follows that
$|A| \leq 2$. We now show that $m$ must be odd.
\medskip

\noindent
{\sc Claim 3:} $m$ is odd. 
\smallskip

Suppose on the contrary that $m$ is even.
Consider the alternating cycle containing the edge $u_0^0 u_1^0$.
It contains vertices $u_2^r, u_3^r, u_4^{r+r3}, \ldots , 
u_{m-1}^{r + r^3 + \cdots + r^{m-3}}, u_0^{r + r^3 + \cdots + r^{m-1} + t}$, etc., where 
$u_0^{r + r^3 + \cdots + r^{m-1} + t}$ is the tail of the two corresponding incident edges on this cycle.
The other alternating cycle containing $u_0^0$ contains vertices
$u_1^1, u_2^1, u_3^{1+r^2}, u_4^{1+r^2},$ $\ldots , u_{m-1}^{1+r^2+\cdots + r^{m-2}}$, 
$u_0^{1+r^2+\cdots + r^{m-2} + t}$, etc., 
where $u_0^{1+r^2+\cdots + r^{m-2} + t}$ is the head of the two corresponding incident edges on this cycle.
Observe that equation (\ref{eq:leI+-rt}) implies that
$r(1+r^2+\cdots + r^{m-2} + t) = r+ r^3 + \cdots + r^{m-1} + t$, and so the vertices
$u_0^{r + r^3 + \cdots + r^{m-1} + t}$ and $u_0^{1+r^2+\cdots + r^{m-2} + t}$ are both contained in $A$. 
Since $|A| = 2$, it follows that $1+r^2+\cdots + r^{m-2} + t = r+ r^3 + \cdots + r^{m-1} + t$ and in addition
either $1+r^2+\cdots + r^{m-2} + t = 0$ or 
$1+r^2+\cdots + r^{m-2} + t = \frac{n}{2}$ with $n$ even. 
But in both cases equation (\ref{eq:leI+-cond}) holds, and so the mapping $\psi$ defined as in (\ref{eq:leI+-psi})
is an automorphism of $X$ which is impossible. Therefore, $m$ is odd, as claimed.
\medskip

\noindent
{\sc Claim 4:} $\mmod{n}{2}{4}$.
\smallskip

Let $C_1$ denote the alternating cycle containing the edge $u_0^0 u_1^0$.
Since $m$ is odd, the vertices  $u_0^{r+r^3+ \cdots + r^{m-2} + t}$
and $u_0^{1+r+r^2+\cdots + r^{m-1} + 2t}$ are contained in $C_1$ with 
$u_0^{r+r^3+ \cdots + r^{m-2} + t}$ being the head of the two corresponding incident edges on $C_1$.
Let $C_2$ denote the other alternating cycle containing $u_0^0$. Then
$u_0^{1+r^2+\cdots + r^{m-1} + t}$ and
$u_0^{1+r+r^2+\cdots + r^{m-1}+2t}$ are vertices of $C_2$ with
$u_0^{1+r^2+\cdots + r^{m-1} + t}$ being the tail of the two corresponding incident edges on $C_2$.
Since $|A| \leq 2$, we thus have that $1+r+r^2+\cdots + r^{m-1} + 2t$ is equal either to
$0$ or to $\frac{n}{2}$, where in the latter case $n$ must be even.
As the former contradicts  half-arc-transitivity of $X$ (see the argument immediately after
equation (\ref{eq:leI+-cond})), we have that $n$ is even and that 
%%%%% equation
\begin{equation}\label{eq:leI+-eq}
1+r+r^2+\cdots + r^{m-1} + 2t = \frac{n}{2}.
\end{equation}
%%%%%
Since $r\in \ZZ_n^*$, $r$ is odd. But this implies that
$1+r+(r^2+r^3)+\cdots  + (r^{m-3} + r^{m-2}) + r^{m-1} + 2t$ is odd too,
and so equation (\ref{eq:leI+-eq}) implies that $\mmod{n}{2}{4}$, as claimed.
\medskip

\noindent
{\sc Claim 5:} $C_0 = u_0^0 u_1^0 u_2^r u_1^r u_0^r u_1^{1+r} u_2^{1+r} u_1^1$ is an $8$-cycle of $X$.
\smallskip

We only need to see that the cardinality of the set $\{0,1,r,r+1\}$ is $4$,
that is, we need to see that $r \neq \pm 1$.
Note first that $r \neq 1$ for otherwise $X$ would be  a Cayley graph of an abelian group and 
thus arc-transitive. Furthermore, $r \neq -1$ for then
$r^m = -1 \neq 1$, contradicting (\ref{eq:leI+-r^m}). 
(Note that $n \neq 2$ for otherwise  $X$ would be isomorphic to 
a lexicographic product of a cycle and $2K_1$, and thus clearly arc-transitive.) 
\medskip

The $8$-cycles belonging to the $\langle \rho, \sigma \rangle$-orbit of $C_0$
will be called the {\em generic} $8$-cycles of $X$. 
We now investigate which 
$8$-cycles, apart from the generic ones, are contained in the $\Aut X$-orbit of $C_0$. 
We assume first that $m \geq 5$, as in this case, since $m$ is odd, 
no $8$-cycle  containing edges from every subgraph $[X_i,X_{i+1}]$, $i \in \ZZ_m$, exists.

By Claim 2 there exists an automorphism $\varphi \in \Aut X$, fixing $u_1^1$
and $u_0^0$ but interchanging $u_2^{1+r}$ and $u_0^1$. We either have $u_1^0 \varphi = u_1^0$ or
$u_1^0\varphi = u_{m-1}^{-r^{m-1}-t}$.
Suppose first that $\varphi$ fixes $u_1^0$. Then it also fixes $u_2^r$. 
Now  since $C_0\varphi$ is an $8$-cycle and since $u_1^{1+r}$ is the tail of both of its incident edges
on $C_0$, we must have  $u_1^{1+r}\varphi = u_1^2$.
Therefore, $u_0^r\varphi = u_2^2$. This leaves us with two possibilities for $u_1^r\varphi$. 
If $u_1^r\varphi = u_1^{2-r}$, then $u_1^{2-r} \to u_2^r$, and
so $2(r-1) = 0$. But $r$ is odd, so that Claim~4 implies $r-1 = 0$, a contradiction. Thus
$u_1^r\varphi = u_3^2$ and so $u_2^r \ot u_3^2$, which forces 
%%%%% equation
\begin{equation}\label{eq:leI+-cond2}
2-r-r^2=0.
\end{equation}
%%%%%
Suppose now that $u_1^0\varphi = u_{m-1}^{-r^{m-1}-t}$. It follows that $u_2^r\varphi = u_{m-2}^{-r^{m-1}-t}$ 
and then a similar argument as above shows that
$u_1^{1+r}\varphi = u_{m-1}^{1-t}$, and so $u_0^r\varphi = u_{m-2}^{1-r^{m-2}-t}$. We now have two possibilities.
Either $u_1^r\varphi = u_{m-1}^{1-r^{m-2}-t}$, which by (\ref{eq:leI+-r^m}) implies that (\ref{eq:leI+-cond2}) holds, or
$u_1^r\varphi = u_{m-3}^{1-r^{m-3}-r^{m-2}-t}$, in which case
$1-r^{m-3}-r^{m-2}-t = -r^{m-1}-t$, that is, 
%%%%% equation
\begin{equation}\label{eq:leI+-cond1}
r^3 + r^2  = r+1 \ \mathrm{or}\ \mathrm{equivalently}\ (r-1)(r+1)^2 = 0. 
\end{equation}
%%%%%
We now show that (\ref{eq:leI+-cond1}) cannot hold. Namely, multiplying by $r^{2(i-1)}$ we get that
$r^{2i+1} + r^{2i} = r+1$ for every $i \in \mathbb{N}$.
Furthermore, by Claim~3 there exists some integer $k$ such that $m = 2k+1$, and so 
%%%%% equation
\begin{equation}\label{eq:leI+-sum}
1 + r + r^2 + \cdots + r^{m-2} + r^{m-1} = k(1+r) + r^{m-1}.
\end{equation}
%%%%%
By (\ref{eq:leI+-r^m}) we therefore have  
$$ 0 = (1+r+\cdots + r^{m-1})(r-1) = k(r^2-1) + r^{m-1}(r-1).$$
Multiplying by $r$ and using (\ref{eq:leI+-cond1}) we obtain
%%%%% equation
\begin{equation}\label{eq:leI+-auxil}
0 = k(r^3-r)+r-1 = k(1-r^2)+r-1 = -k(r-1)(r+1)+(r-1).
\end{equation}
%%%%%
Now, since $r-1 \neq 0$ is even, Claim~4 and equation (\ref{eq:leI+-cond1}) imply
that there exists some odd prime $q$ dividing $r+1$ and $n$ but not $r-1$.
However, equation (\ref{eq:leI+-auxil}) implies that $q$ does divide $r-1$, a contradiction.
Therefore, (\ref{eq:leI+-cond1}) cannot hold, and so
(\ref{eq:leI+-cond2}) holds.

Again using Claim 2, there also exists an automorphism $\vartheta \in \Aut X$ fixing $u_0^0$ and $u_1^0$ but
interchanging $u_1^1$ and $u_{m-1}^{-t}$. An  analysis similar to the one used above 
shows that the only possibilities for the image  $C_0\vartheta$ are:
$$
\begin{array}{ll}
u_0^0 u_1^0 u_2^r u_3^{r+r^2} u_2^{r+r^2} u_1^{r+r^2} u_0^{-1+r+r^2} u_{m-1}^{-t}& ,\\
u_0^0 u_1^0 u_2^r u_1^{r} u_2^{2r} u_1^{2r} u_0^{-1+2r} u_{m-1}^{-t} & \mathrm{and}\\
u_0^0 u_1^0 u_2^r u_1^{r} u_0^{r} u_{m-1}^{r-t} u_{m-2}^{r-r^{m-2}-t} u_{m-1}^{-t}. &
\end{array}
$$
The conditions for the above $8$-cycles to exist are, respectively,
(\ref{eq:leI+-cond1}), 
%%%%% equation
\begin{equation}\label{eq:leI+-cond2'}
1+r-2r^2 = 0
\end{equation}
%%%%%
and 
%%%%% equation
\begin{equation}\label{eq:leI+-cond3}
r^3 = 1.
\end{equation}
%%%%%
Recall that (\ref{eq:leI+-cond1}) cannot hold and that (\ref{eq:leI+-cond2}) does hold. 
It is easy to check that (\ref{eq:leI+-cond2}) and (\ref{eq:leI+-cond2'}) imply
(\ref{eq:leI+-cond3}) and that (\ref{eq:leI+-cond2}) and (\ref{eq:leI+-cond3}) imply 
(\ref{eq:leI+-cond2'}). Thus all three conditions hold. From (\ref{eq:leI+-cond2})
and (\ref{eq:leI+-cond2'}) we get that
%%%%% equation
\begin{equation} 
\label{eq:leI+-r-1^2}
(r-1)^2 = 0.
\end{equation}
%%%%%
Then 
$C = u_0^0 u_1^0 u_2^r u_3^r u_2^{r-r^2} u_1^{r-r^2} u_2^{2r-r^2} u_1^{2r-r^2} =
u_0^0 u_1^0 u_2^r u_3^r u_2^{1-r} u_1^{1-r} u_2^{1} u_1^{1} $ is 
an $8$-cycle of $X$ (recall that $r \neq 1$). By Claim 2 there exists an automorphism
$\eta \in \Aut X$ fixing $u_2^r$ and $u_3^r$ but interchanging $u_2^{r-r^2}$ with $u_4^r$. 
But this implies $u_1^0\eta = u_1^0$, $u_0^0\eta = u_0^0$ and $u_1^{r-r^2} \eta = u_5^{r+r^4}$, and so
$C\eta$ is not an $8$-cycle, a contradiction. 
This shows that $X$ cannot exists when $m \geq 5$.

Since $m$ is odd this leaves us with  $m=3$ as the only  other possibility.
A similar analysis as in the general case shows that the only possibilities for $C_0\varphi$
(where $\varphi \in \Aut X$ fixes $u_1^1$ and $u_0^0$ but interchanges $u_2^{1+r}$ and $u_0^1$) are 
$8$-cycles which exist only when (\ref{eq:leI+-cond2}) holds. (In this analysis we get that the
only possibilities not encountered in the general case are those for which $2r+r^2+2t=0$ or $2+r^2+2t=0$,
which are of course both impossible as $n$ is even.)
By (\ref{eq:leI+-r^m}) we have $r^3 = 1$, and so multiplying by $r$ in (\ref{eq:leI+-cond2}) we get
that (\ref{eq:leI+-r-1^2}) holds.
Thus the $8$-cycle $C = u_0^0 u_1^0 u_2^r u_0^{r+t} u_2^{1-r} u_1^{1-r} u_2^{1} u_1^{1}$ exists in $X$.
Again let $\eta \in \Aut X$ be an automorphism fixing $u_2^r$ and $u_0^{r+t}$ but interchanging 
$u_2^{1-r}$ with $u_1^{r+t}$. Then $u_1^0\eta = u_1^0$, $u_0^0\eta = u_0^0$ and $u_1^{1-r}\eta = u_2^{2r+t}$.
Note that we cannot have $u_1^1\eta = u_1^1$ for otherwise $u_2^1\eta = u_2^1$, which contradicts
the fact that $u_2^1 \sim u_1^{r-r^2}$. It follows that $u_1^1\eta = u_2^{-t}$, and so $u_2^1\eta = u_1^{-r-t}$.
Consequently $-r-t = 2r+t$, that is $3r+2t=0$, which forces $1+r+r^2+2t=0$.
As this contradicts (\ref{eq:leI+-eq}), the proof is complete.
\end{proof}
\bigskip

We are now ready to prove Theorem~\ref{the:TA=2}.
\bigskip

\noindent
{\sc Proof of Theorem~\ref{the:TA=2}:}\\
That connected quartic half-arc-transitive weak metacirculants of Class~I are tightly attached now
follows by Lemmas~\ref{le:class1++} and \ref{le:class1+-}. To prove the converse observe,
as already noted in Section~\ref{sec:examples}, that the results of \cite{DM98,PSxx} imply that every
connected quartic tightly attached half-arc-transitive graph is isomorphic either to some
$\X_o(m,n;r)$ or to some $\X_e(m,n;r,t)$. As these two graphs are clearly both weak metacirculants of Class~I,
the proof is complete.
\hfill $\Qed$

%%%%%%%%%%%%%%%%%%%%%%%%%%%%%%%%%%%%%%%%%%%
%%%%%%%%%%%%%%%%%%%%%%%%%%%%%%%%%%%%%%%%%%%
%%%%%%%%%%%%%%%   Class II %%%%%%%%%%%%%%%%
%%%%%%%%%%%%%%%%%%%%%%%%%%%%%%%%%%%%%%%%%%%
%%%%%%%%%%%%%%%%%%%%%%%%%%%%%%%%%%%%%%%%%%%

\section{Graphs of Class~II}
\label{sec:ClassII}

\indent

In this section connected quartic
half-arc-transitive metacirculants of Class~II are studied in great detail.
The following theorem is the main result.

%%%%%  theorem the:IItheorem  %%%%%
\begin{theorem}
\label{the:IItheorem}
Let $X$ be a connected quartic half-arc-transitive weak ($m$,$n$)-metacirculant of Class~II. 
Then the following hold:
\begin{enumerate}
\item[(i)] $X$ is a Cayley graph for the group $\la \rho, \sigma\ra$, where $(\rho,\sigma)$ is some pair of automorphisms
of $X$ such that $X$ is a weak $(m,n)$-metacirculant of Class~II relative to $(\rho,\sigma)$,
\item[(ii)] $(\Aut X)_v \iso \ZZ_2$ for all $v \in V(X)$,
\item[(iii)] $m$ divides $n$ and $d_m  = \frac{n}{m} > 2$,
\item[(iv)] there exist $r \in \ZZ_n^*$ and $t \in \ZZ_n$ such that $X \cong \Y(m,n;r,t)$, where
	parameters $r$ and $t$ satisfy the following conditions:\\
	$\bullet$ $r^m = 1$,\\
	$\bullet$ $m(r-1) = t(r-1) = (r-1)^2 = 0$,\\
	$\bullet$ $\la m \ra = \la t \ra$ in $\ZZ_n$,\\
	$\bullet$ there exists a unique $c \in \{0,1,\ldots , d_m - 1\}$ such that $t = cm$ and $m = ct$ and\\
	$\bullet$ there exists a unique $a \in \{0,1,\ldots , d_m - 1\}$ such that $at = -am = r-1$,
\item[(v)] $X$ is tightly attached unless $m$ and $d_m$ are both even, 
and $n = 8n_1$ for some integer $n_1 > 2$, where $n_1$ is even or not squarefree.
\end{enumerate}
\end{theorem}

Then, using this result, we present a list of all connected quartic half-arc-transitive 
weak metacirculants of Class~II
of order up to $1000$ which are not tightly attached  (see Table~\ref{table:IInotta}). Finally, 
we construct an infinite family of such graphs (see Construction~\ref{cons:infinite}).

Throughout this section we let $X$ denote a connected quartic half-arc-tran\-sitive
weak $(m,n)$-metacirculant of Class~II. Choose some automorphisms $\rho$ and $\sigma$ such that 
$X$ is a weak $(m,n)$-metacirculant of Class~II relative to the ordered pair $(\rho,\sigma)$.
Fix a vertex $u \in V(X)$ and let
$u_i^0 = u\sigma^i$ for all $i \in \{0,1,\ldots , m-1\}$. 
Then let
$u_i^j = u_i^0\rho^j$ for all $i \in \ZZ_m$, $j \in \ZZ_n$.
Thus $X_i = \{u_i^j\ |\ j \in \ZZ_n \}$,
$i \in \ZZ_m$, are the orbits of $\rho$ and $X_i = X_0\sigma^i$.
We shall say that an edge connecting vertices from the same orbit
$X_i$ is an {\sl inner edge} and that an edge connecting vertices from different
orbits is an {\sl outer edge}.

Since $d_{inn}(X) = 2$,
there exists some nonzero $s \in \ZZ_n$ such that
$u_0^j \sim u_0^{j \pm s}$ for all $j \in \ZZ_n$.
Fix an orientation of edges induced on $X$ by the half-arc-transitive
action of $\Aut X$ and denote the corresponding directed graph by $D_X$.
Then the indegrees and the outdegrees of the subgraphs of $D_X$ 
induced  by $X_i$ are all equal to $1$.
We will assume that
$u_0^{j-s} \rightarrow u_0^j \rightarrow u_0^{j+s}$.
Letting $r \in \ZZ_n^*$ be as in equation (\ref{eq:metagroup}), we have that
$u_0^s\sigma^i = u_0^0\rho^s\sigma^i = u_0^0\sigma^i\rho^{r^i s}$, and so
\begin{equation}\label{eq:sI}
u_i^{j - r^i s} \rightarrow u_i^j \to u_i^{j + r^i s},
\quad \mathrm{for}\ \mathrm{all}\ i \in \ZZ_m,\ j \in \ZZ_n.
\end{equation}

There exists some $k \in \ZZ_m \setmin \{0\}$ such that
the vertices from the orbit $X_0$ are adjacent to the vertices from the
orbit $X_k$. Since $X$ is connected, $\la k \ra = \ZZ_m$, so that we can
assume $k = 1$ (otherwise
take $\sigma' = \sigma^k$ and $r' = r^k$). Let $a \in \ZZ_n$ be such that
$u_0^0 \sim u_1^a$. We can assume that $u_0^0 \to u_1^a$
(otherwise take $\rho' = \rho^{-1}$ and then
choose the other of the two possible orientations of the edges for $D_X$).
With no loss of generality we can also assume that $a = 0$ 
(otherwise take $\sigma' = \sigma\rho^a$). 
Therefore, $u_i^j \to u_{i+1}^{j}$ 
for all $i \in \ZZ_m \setmin \{m-1\}$, $j \in \ZZ_n$.
Since $\sigma$ cyclically permutes the $m$ orbits of $\rho$, we have that
$u_0^0\sigma^m \in X_0$.
Thus, there exists a unique $t \in \ZZ_n$ such that
$u_0^0\sigma^m\rho^{-t} = u_0^0$. Since
the orbits $X_i$ are blocks of imprimitivity for the group
$H = \langle \rho, \sigma \rangle$, half-arc-transitivity of $X$ and the 
orientation of the edges of $D_X$ imply that 
an element of $H$ fixing a vertex must necessarily
fix all of its neighbors pointwise. By connectedness of $X$ we then have that
%%%%% equation %%%%%
\begin{equation}
\label{eq:IIregularly}
\mathrm{the\ group\ }\la \rho, \sigma\ra \mathrm{\ acts\ regularly\ on\ }V(X).
\end{equation}
%%%%%
In particular, $\sigma^m = \rho^t$. This implies that
$\rho = \rho^{-t}\rho\rho^t = \sigma^{-m}\rho\sigma^m = \rho^{r^m}$, and so
%%%%% equation %%%%%
\begin{equation}\label{eq:IIr^m}
r^m = 1.
\end{equation}
%%%%%
Moreover, $u_i^j\sigma = u_i^0\rho^j\sigma = u_i^0\sigma\rho^{rj} =
u_{i+1}^0\rho^{rj} = u_{i+1}^{rj}$ for $i \neq m-1$, $j \in \ZZ_n$. By (\ref{eq:IIr^m}) we 
now also have
$u_{m-1}^j\sigma = u_0^{rj}\sigma^{m-1}\sigma = u_0^{rj}\rho^t = u_0^{rj+t}$, and so
%%%%% equation %%%%%
\begin{equation}
\label{eq:IIsigma}
u_i^j\sigma = \left\{\begin{array}{ccl} u_{i+1}^{rj} & ; & i \in \ZZ_m \setmin \{m-1\},\ j \in \ZZ_n\\
				    u_0^{rj+t} & ; & i = m-1,\ j \in \ZZ_n \end{array}\right.\quad \quad\mathrm{and} 
\end{equation}
%%%%%
%%%%% equation %%%%%
\begin{equation} \label{eq:IIsigmaarcs}
u_i^j \to \left\{\begin{array}{ccl} u_{i+1}^{j} & ; & i \in \ZZ_m \setmin \{m-1\},\ j \in \ZZ_n \\
				    u_0^{j+t} & ; & i = m-1,\ j \in \ZZ_n .\end{array}
				    \right.  \end{equation}
Let us now consider the edge $u_{m-1}^0u_0^t$. By (\ref{eq:IIsigma}), $\sigma$
maps this edge to the edge $u_0^t u_1^{rt}$, and so (\ref{eq:IIsigmaarcs}) implies that $rt = t$, that is,
%%%%% equation %%%%%
\begin{equation}\label{eq:IIrt=t}
t(r-1) = 0.
\end{equation}
%%%%%

We claim that $rs \neq \pm s$.
Suppose on the contrary that $rs = s$ or $rs = -s$ and consider the
permutation $\varphi$ of $V(X)$ defined by 
the rule $u_0^j\varphi = u_0^{-j}$, where $j \in \ZZ_n$, and
$u_i^j\varphi = u_{m-i}^{-j-t}$, where
$i \in \ZZ_m \setmin \{0\}$ and $j \in \ZZ_n$. Since $rs = \pm s$,
we  have that either $r^is = s$ or that $r^is = (-1)^{i}s$ for all $i \in \ZZ_m$. It is now easy 
to check that $\varphi$ is an automorphism of $X$. 
But $\varphi\rho^s$ interchanges adjacent vertices $u_0^0$ and $u_0^s$, which by Proposition~\ref{pro:flip}
contradicts half-arc-transitivity of $X$.

We now investigate certain $8$-cycles of $X$ in order to obtain a better 
understanding of the structural properties of $X$.
Consider the following closed walk of $X$:
%% the generic 8-cycle
\begin{equation} \label{eq:II8cycle}
(u_0^0, u_0^s, u_1^s, u_1^{s+rs}, u_0^{s+rs}, u_0^{rs}, u_1^{rs}, u_1^0, u_0^0).
\end{equation}
Since  $s \neq 0$, $r \in \ZZ_n^*$ and $rs \neq \pm s$, it follows that the above $8$
vertices are all distinct, and so the closed walk~(\ref{eq:II8cycle}) gives rise to an $8$-cycle. 
Every $8$-cycle of $X$ belonging to the $H$-orbit of this $8$-cycle will be called a
{\em generic} $8$-cycle. 

\begin{figure}
\begin{center}
\includegraphics[scale=0.4]{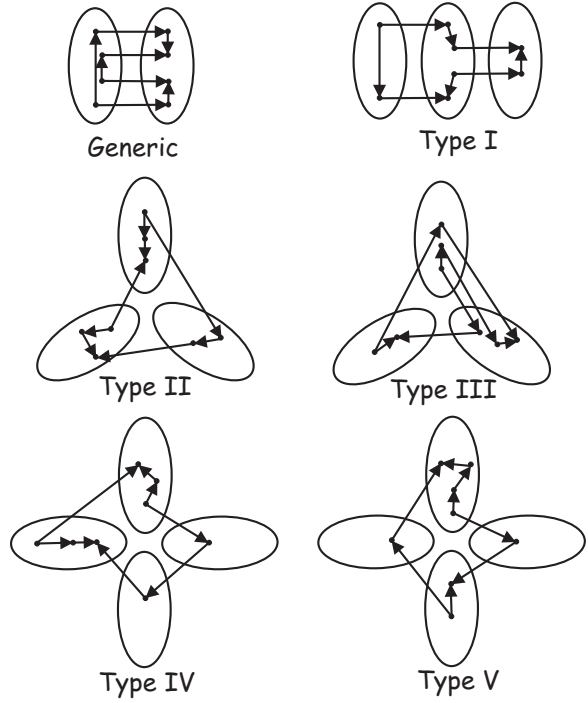}
\caption{Possible types of $8$-cycles having code $11100100$.}
\label{fig:8cycles}
\end{center}
\end{figure}

To every $8$-cycle $C$ of $X$ we
assign a binary sequence as follows. When traversing
$C$, we assign  value $1$ to each edge of $X$ traversed along its orientation in $D_X$,
and we assign value $0$ to each edge of $X$ traversed against its orientation in $D_X$.
We say that two binary sequences corresponding to $8$-cycles of $X$ 
are {\em equivalent} if one can be obtained from the other using cyclic rotations and
reflections. We let the {\em code} of $C$ be the equivalence class of its sequences and we
denote it by any of the corresponding sequences.
Therefore, the code of the generic $8$-cycle given in (\ref{eq:II8cycle}) is $11100100$ 
(see Figure~\ref{fig:8cycles}). Note that,
since $X$ is half-arc-transitive,
the code of a cycle is invariant under the action of $\Aut X$.
 On the other hand, 
there exists an automorphism $\tau \in \Aut X$ fixing
$u_1^{rs}$ and interchanging $u_1^0$ and $u_0^{rs}$. 
Since $u_0^{rs} \to u_0^{s+rs}$, we thus have that $u_0^{s+rs}\tau = u_2^0$. 
Consequently, the image under $\tau$ of the  generic  $8$-cycle corresponding to  (\ref{eq:II8cycle})  is 
an $8$-cycle consisting of vertices from at least
three orbits $X_i$ and is therefore not generic. 
The following lemma gives all possible $H$-orbits of $8$-cycles of $X$ having code $11100100$.

%%%%%  lemma le:IIpossHorbits  %%%%%
\begin{lemma}
\label{le:IIpossHorbits}
	With the notation introduced in this section
	the only possible $H$-orbits of $8$-cycles  having code $11100100$ in $X$ 
	are given in Table~\ref{tab:classIItypes}, together with the corresponding representatives 
	and the necessary and sufficient arithmetic conditions for their existence.
	\end{lemma}

\begin{proof}
Let $C = c_0c_1c_2c_3c_4c_5c_6c_7$ be an $8$-cycle with code $11100100$. We divide our investigation
into several cases depending on the number of orbits $X_i$ the $8$-cycle $C$ meets. 

\medskip
\noindent
{\sc Case 1:} $C$ meets one orbit.
\smallskip

As any such $8$-cycle has code $11111111$, this case cannot occur.
\medskip

\noindent
{\sc Case 2:} $C$ meets two orbits.
\smallskip

Clearly, the number of outer edges of $C$ is even. In fact, $C$ either has $2$ or $4$ outer edges.
The former case is impossible for otherwise $C$ contains at least four consecutive vertices in 
a single orbit and thus $1111$ is a subsequence of the code of $C$.
It is thus clear that the inner and outer edges alternate on $C$. Therefore, the first and the last $1$ of 
the subsequence $111$ of the code of $C$ both correspond to inner edges, and so 
it is clear that $C$ is a generic $8$-cycle.
\medskip

\noindent
{\sc Case 3:} $C$ meets three orbits, say, with no loss of generality, $X_0, X_1$ and $X_2$.
\smallskip

%%% m > 3
Suppose first that $m > 3$. Therefore, if $c_i \in X_0$ or $c_i \in X_2$, at least
one of
$c_{i-1}, c_{i+1}$ lies in the same orbit as $c_i$. 
This implies that no four consecutive vertices of $C$ are contained in a single orbit. 
Namely, they cannot be contained in
$X_0$ or $X_2$, for otherwise the code of $C$ would contain  $1111$ as a subsequence.
Moreover, they cannot be contained in $X_1$ since
there are no edges between $X_0$ and $X_2$. We now show that no three
consecutive vertices of $C$ are contained in a single orbit. 
Suppose on the contrary that $c_0, c_1$ and $c_2$ are all contained in one orbit.
If this orbit is $X_0$, then $c_3, c_7 \in X_1$, so in order to have the required code, 
at least one of $c_4$ and $c_6$ lies in $X_1$. But then the remaining two vertices
lie in $X_2$, so the code cannot be $11100100$. A similar argument shows that the orbit
containing $c_0, c_1$ and $c_2$ cannot be $X_2$. 
Suppose now that the orbit containing $c_0,c_1$ and $c_2$ is $X_1$. 
It is then clear that one of $c_3$ and $c_7$ lies in $X_0$ and the other in $X_2$;
say $c_3 \in X_0$ and $c_7 \in X_2$.
It follows that $c_4 \in X_0$, $c_5 \in X_1$ and $c_6 \in X_2$. It is easy to see however, that
such an $8$-cycle does not have code $11100100$.
Therefore, no three consecutive vertices of $C$ lie on a single orbit.
Consequently, $X_0$ and $X_2$ each contain two vertices of $C$ and so
four vertices of $C$ are contained in $X_1$.
%Without loss of generality assume
%$c_0, c_1 \in X_0$,  $c_2,c_3, c_6,c_7 \in X_1$ and $c_4,c_5 \in X_2$.
As $C$ has code $11100100$, it is now clear that $C$ lies in the $H$-orbit
of the $8$-cycle from row~2 of Table~\ref{tab:classIItypes}.
We say that the $8$-cycles of this $H$-orbit are of {\em type~I} (see Figure~\ref{fig:8cycles}).

%%% m = 3
Suppose now that $m = 3$. With no loss of generality we can assume that the
sequence $11100100$ is obtained when traversing $C$ according to increasing subscripts of vertices
and, in addition,  that $c_0 \in X_0$ and that the walk 
$(c_0,c_1,c_2,c_3)$ gives rise to the subsequence $010$. 
We first consider the possibility that $c_1 \in X_0$.
Then $c_2, c_3 \in X_1$. We claim that this forces $c_4 \in X_1$. 
Namely, if this is not the case,
then $c_4, c_5 \in X_0$, and so the fact that
$c_5 \to c_6 \to c_7 \ot c_0$ implies that $C$ does not contain vertices from $X_2$, a contradiction.
Therefore, $c_4 \in X_1$ and hence $c_5 \in X_2$. It is now easy to see
that the only way for $C$ to have the required code is to have
$c_6 \in X_2$ and $c_7 \in X_0$. 
Thus, $C$ is contained in the $H$-orbit of $8$-cycles
whose representative is given in row~3 of Table~\ref{tab:classIItypes}.
We say that the $8$-cycles of this $H$-orbit are of {\em type~II} (see Figure~\ref{fig:8cycles}).
Consider now the possibility that $c_1 \notin X_0$, and so $c_1 \in X_2$. 
It follows that  $c_2 \in X_2$ and $c_3 \in X_1$. 
If $c_4 \in X_1$, then 
$c_5 \in X_2$, and then the only way for $C$ to have the required code is to have
$c_6 \in X_2$ and $c_7 \in X_0$. Note
that this $8$-cycle is of type~I.
If however $c_4 \in X_0$, then $c_5 \in X_0$, and then the only way for $C$ to have the required code is 
to have $c_6, c_7 \in X_1$.
Thus $C$ is in the $H$-orbit of the $8$-cycle from row~4 of Table~\ref{tab:classIItypes}.
We say that such  $8$-cycles are of {\em type~III} (see Figure~\ref{fig:8cycles}). 
To summarize, if $m = 3$ we can have up to four different types of $8$-cycles
with code $11100100$: the generic ones and $8$-cycles of types I, II and III.

\medskip
\noindent
{\sc Case 4:} $C$ meets four orbits,
say, with no loss of generality, $X_0, X_1, X_2$ and $X_3$.

\smallskip

Observe first that no $8$-cycle with code $11100100$ exists if $m > 4$. 
Namely, in this case an  $8$-cycle would
necessarily have to contain precisely two vertices from each of the orbits $X_0$, $X_1$, $X_2$ and $X_3$.
But then the code of $C$ would contain $1111$ as a subsequence, which is impossible.
We can therefore  assume that $m = 4$.
Because of the particular code of $C$ it is clear that $C$ has at most three consecutive outer edges.
Hence, either all the outer edges of $C$ give rise to digit $1$ in  the code of $C$ or
they all give rise to digit $0$.
It follows that  $C$ has precisely four inner and four outer edges.
If the outer edges give rise to digit $1$ in  the code  $11100100$ of $C$, 
then $C$ belongs to the $H$-orbit of the
$8$-cycle from row~5 of Table~\ref{tab:classIItypes}. 
We say that such  $8$-cycles are of {\em type~IV} (see Figure~\ref{fig:8cycles}).
If on the other hand the outer edges 
give rise to digit $0$ in  the code  $11100100$ of $C$, 
then $C$ belongs to the $H$-orbit of the
$8$-cycle from row~6 of Table~\ref{tab:classIItypes}. 
We say that the $8$-cycles of this
$H$-orbit are of {\em type~V} (see Figure~\ref{fig:8cycles}).

\medskip
\noindent
{\sc Case 5:} $C$ meets more than four orbits.

\smallskip

It is easy to see that no such $8$-cycle exists.
\end{proof}

\begin{table}[!htb]
\begin{footnotesize}
\begin{center}
\begin{tabular}{@{}c|l|l|c@{}}
	Row	&  Type   &  A representative 	&  Condition	\\
	\hline  & & & \\
   $1$ & generic  & $u_0^0 u_0^s u_1^s u_1^{s+rs} u_0^{s+rs} u_0^{rs}u_1^{rs} u_1^0$ & none \\    
   $2$ & type I   & $u_0^0u_1^0u_1^{sr}u_2^{sr}u_2^{s(r-r^2)}u_1^{s(r-r^2)}u_1^{s(2r-r^2)}u_0^{s(2r-r^2)}$ & $s(1-2r+r^2) = 0$ \\
   $3$ & type II  & $u_0^0u_1^{0}u_1^{sr}u_2^{sr}u_2^{s(r-r^2)}u_2^{s(r-2r^2)}u_0^{s(r-2r^2)+t}u_0^{s(-1+r-2r^2)+t}$ & $s(2-r+r^2)-t=0$ and $m = 3$ \\
   $4$ & type III & $u_0^0u_0^{s}u_1^su_1^{s(1+r)}u_0^{s(1+r)}u_2^{s(1+r)-t}u_2^{s(1+r+r^2)-t}u_1^{s(1+r+r^2)-t}$ & $s(1+r+r^2)-t=0$ and $m = 3$ \\
   $5$ & type IV  & $u_0^0u_1^0u_2^0u_3^0u_3^{-sr^3}u_3^{-2sr^3}u_0^{-2sr^3+t}u_0^{s(-1-2r^3)+t}$ & $s(2+2r^3)-t=0$ and $m = 4$ \\
   $6$ & type V   & $u_0^0u_0^su_0^{2s}u_0^{3s}u_3^{3s-t}u_2^{3s-t}u_2^{s(3+r^2)-t}u_1^{s(3+r^2)-t}$ & $s(3+r^2)-t=0$ and $m = 4$ \\   
\end{tabular}
\caption{Possible $H$-orbits of $8$-cycles of code $11100100$.}
\label{tab:classIItypes}
\end{center}
\end{footnotesize}
\end{table}

%%%%%  proposition pro:stab  %%%%%
\begin{proposition}\label{pro:IIstab}
Let $X$ be a connected quartic half-arc-transitive weak ($m$,$n$)-meta\-circulant of Class~II and 
let $v \in V(X)$. Then $(\Aut X)_v \cong \ZZ_2$ or possibly  $(\Aut X)_v \cong \ZZ_2 \times \ZZ_2$ in which case
$m = 4$.
\end{proposition}

\begin{proof}
Let $\rho, \sigma \in \Aut X$ be such that $X$ is a weak $(m,n)$-metacirculant of Clas~II relative to 
the ordered pair $(\rho, \sigma)$ and that all the assumptions made in the third paragraph of this section
hold. Moreover, adopt the notation introduced in this section and 
let $C$ denote the generic $8$-cycle 
from row~1 of Table~\ref{tab:classIItypes}. We distinguish two cases
depending on whether $m$ equals $4$ or not.

\medskip

\noindent
%%% m \neq 4                         
{\sc Case 1:} $m \neq 4$.

\smallskip

If the stabilizer $(\Aut X)_{u_0^s}$ is not isomorphic to $\ZZ_2$, 
then there exists an automorphism $\varphi$ of $X$ which fixes $u_0^0$ and $u_0^s$, and maps
$u_1^s$ to $u_0^{2s}$. Therefore, $C\varphi$ is an $8$-cycle (with code $11100100$) containing
three consecutive vertices ($u_0^0, u_0^s$ and $u_0^{2s}$)
in a single orbit of $\rhogr$. It follows, by Lemma~\ref{le:IIpossHorbits}, that
$m = 3$ and that $C\varphi$ is of type~II. Thus we must
have $u_1^{s+rs}\varphi = u_2^{2s-t}$. 
But this is impossible since then $u_1^s \to u_1^{s+rs}$ and $u_1^s\varphi \ot u_1^{s+rs}\varphi$.
It follows that $(\Aut X)_v \cong \ZZ_2$, as claimed.

\medskip

\noindent
%%% m = 4                            
{\sc Case 2:} $m = 4$.

\smallskip

Suppose that $|(\Aut X)_v| > 2$. 
Then there exists an automorphism
$\varphi$ of $X$ fixing $u_0^{s+rs}$ and $u_1^{s+rs}$, and interchanging
$u_0^{rs}$ and $u_3^{s+rs-t}$. This implies that $u_1^{rs}\varphi = u_3^{s+rs+r^3s-t}$, 
$u_1^0\varphi = u_2^{s+rs+r^3s-t}$ and $u_1^s\varphi = u_1^s$. It follows that
$C\varphi$ is of type~V. Therefore, $8$-cycles of type~V exist in $X$. 

To complete the proof we now show that the only automorphism of
$X$ fixing a vertex and all of its neighbors is the identity. 
To this end let $\varphi \in \Aut X$ be an automorphism fixing $u_0^0$ and its four neighbors
$u_0^{-s}, u_0^s, u_{3}^{-t}$ and $u_1^0$. 
There exists a unique $8$-cycle $C'$ with code $11100100$ containing vertices $u_0^{-s},u_0^0,u_0^s$ 
and $u_0^{2s}$. (It is of type~V.)
Since $\varphi$ fixes the first three of these four vertices, 
and since $u_0^{-s} \to u_0^0 \to u_0^s \to u_0^{2s}\varphi$ is
a directed path of $C'\varphi$, the $8$-cycle $C'\varphi$ is of type~V. 
Consequently,
$\varphi$ fixes all of its vertices pointwise. In particular $u_0^{2s}\varphi = u_0^{2s}$.
It is now clear that $\varphi$ fixes $u_0^s$ and all 
of its neighbors. Continuing inductively, we see that 
$\varphi$ fixes every vertex of form $u_0^{js}$ and all of its neighbors.
Considering again  the generic $8$-cycle $C$.
Since $\varphi$ fixes its vertices $u_0^0, u_0^s, u_0^{rs}, u_0^{s+rs}$ and all
of their neighbors, it fixes $C$ pointwise. It follows that $\varphi$
fixes $u_1^0$ and all of its neighbors. Since $X$ is connected a repeated use
of the above argument finally shows that $\varphi$ is the identity, as required. 
It is now clear, that $|(\Aut X)_v|= 4$, and so $\Aut X \cong \ZZ_2 \times \ZZ_2$, as claimed.
\end{proof}\bigskip

\noindent
{\bf Remark:} In fact, as we shall see in Theorem~\ref{the:IItheorem}, the vertex stabilizer
cannot be isomorphic to $\ZZ_2 \times \ZZ_2$.
\bigskip

%%%%%  proposition pro:IImdivn  %%%%%
\begin{proposition}\label{pro:IImdivn}
Let $X$ be a connected quartic half-arc-transitive weak ($m$,$n$)-meta\-circulant of Class~II.
Then $m$ divides $n$ and moreover, 
there exist $r \in \ZZ_n^*$ and $t \in \ZZ_n$, satisfying (\ref{eq:IIr^m}) and (\ref{eq:IIrt=t}), 
such that $X \cong \Y(m,n;r,t)$.
\end{proposition}

\begin{proof}
Let $\rho, \sigma \in \Aut X$ be such that $X$ is a weak $(m,n)$-metacirculant of Clas~II relative to 
the ordered pair $(\rho, \sigma)$ and that all the assumptions made in the third paragraph of this section
hold. Moreover, adopt the notation introduced in this section.
Let $d_s$ denote the
order of $s$ in $\ZZ_n$. There exist unique integers 
$a \geq 0$ and $b \in \{0,1,\ldots ,m-1\}$ such that 
$d_s = am + b$.
Let $C_0$ denote the directed $d_s$-cycle 
$u_0^0 u_0^s u_0^{2s} \cdots u_0^{(d_s-1)s}$. 
By Proposition~\ref{pro:IIstab} there exists a unique automorphism $\tau \in \Aut X$, which
fixes $u_0^0$, interchanges $u_0^s$ and $u_1^0$, and interchanges $u_0^{-s}$ and $u_{m-1}^{-t}$.
We claim, that the
image $C_0\tau$ of $C_0$ under $\tau$ is the directed $d_s$-cycle at $u_0^0$ consisting
only of outer edges. 
Suppose this does not hold.
Then there exists a smallest $k \in \{1,2, \ldots ,d_s-1\}$ such that 
$\tau$ maps the inner edge $u_0^{ks} u_0^{(k+1)s}$ to an inner edge. 
Since $H = \la \rho, \sigma \ra$ acts transitively on $V(X)$, there exists an automorphism
$\varphi \in H$ such that $u_0^{ks}\tau = u_0^{ks}\varphi$. 
The orbits $X_i$ of $\rho$ are blocks of imprimitivity for $H$, and so it is clear
that $\varphi$ maps inner edges to inner edges. Therefore, we also have
$u_0^{(k+1)s}\tau = u_0^{(k+1)s}\varphi$. However, as $\tau \notin H$, 
$\tau\varphi^{-1}$ is a nontrivial automorphism of $X$ fixing an edge.
Hence, Proposition~\ref{pro:IIstab} implies that $(\Aut X)_{u_0^0} \iso \ZZ_2 \times \ZZ_2$ and that $m = 4$.
Moreover, following its proof we see  that $8$-cycles with code $11100100$ of type V exist. 
In particular there exists a unique $8$-cycle $C_1$ of type~V containing
vertices $u_0^{-s}, u_0^0,u_0^s$ and $u_0^{2s}$. Since $\tau$ maps the first three 
vertices to $u_3^{-t}, u_0^0$ and $u_1^0$, respectively, it is clear that $C_1\tau$ is of type~IV.
It follows that $u_0^{2s}\tau = u_2^0$. We now repeatedly use this argument 
on $8$-cycles of type~V containing vertices $u_0^{(i-1)s}, u_0^{is}, u_0^{(i+1)s}$
and $u_0^{(i+2)s}$ to finally prove that the edge $u_0^{ks} u_0^{(k+1)s}$ gets mapped 
to an outer edge, a contradiction which proves our claim.

Observe that the fact that $C_0\tau$ is the directed $d_s$-cycle at $u_0^0$ consisting
only of outer edges implies that $u_b^{at} = u_0^0$, and so $b = 0$ and $at = 0$.
In particular, this shows that the order of $t$ in $\ZZ_n$ divides $a$.
Since $d_s = am$ is the order of $s$ in $\ZZ_n$, we have that the order of $t$ in $\ZZ_n$ divides $d_s$. 
Therefore, $\langle t \rangle$ is a subgroup of $\langle s \rangle$. 
However, the connectedness of $X$ implies that 
$\langle s, t \rangle = \ZZ_n$, and so $\langle s \rangle = \ZZ_n$, that is $\gcd(n,s) = 1$.
It is now clear that $X \cong \Y(m,n;r,t)$. Finally, 
the equation $n = d_s = am$ implies that $m$ divides $n$.
\end{proof}\bigskip

\noindent
{\bf Remark:} Note that Proposition~\ref{pro:IImdivn} implies that
one can assume $s = 1$ in (\ref{eq:sI}), that is, $u_0^0 \to u_0^1$.
For the rest of this section we therefore let $s = 1$.
\bigskip

For future reference we record the nature of the action of the automorphism $\tau$ from the
proof of the above proposition.

%%%%%  lemma le:IItau  %%%%%
\begin{lemma}\label{le:IItau}
With the notation introduced in this section let $d_m$ be the unique 
integer such that $n = md_m$.
Then for every $i \in \ZZ_m$ and every $j \in \ZZ_n$ there exist unique
integers $a \in \{0,1,\ldots ,d_m - 1\}$ and $b \in \{0,1,\ldots , m-1\}$
such that $j = (am + b)r^i$ in $\ZZ_n$. Moreover, the 
unique automorphism $\tau$ of $X$ fixing $u_0^0$,
interchanging $u_0^1$ and $u_1^0$, and interchanging $u_0^{-1}$ and $u_{m-1}^{-t}$,
maps according to the rule $u_i^j \tau = u_b^{i+at}$.
%, where $j = (am + b)r^i$ is as above.
\end{lemma}

\begin{proof}
Observe first that since $r \in \ZZ_n^*$ the existence of unique $a$
and $b$ is clear. The proof of Proposition~\ref{pro:IImdivn}
shows that $\tau$ maps the inner edges of $X_0$ to outer edges. 
Therefore it maps the outer edges connecting $X_0$ to $X_1$ to inner edges. 
Continuing inductively we can  see that $\tau$ interchanges
inner edges with outer edges. It is now clear that $u_i^0\tau = u_0^i$ and that $u_i^{r^i}\tau = u_1^i$, 
$u_i^{2r^i}\tau = u_2^i$, etc. Finally, $u_i^{(am + b)r^i}\tau = u_b^{i+at}$, which completes the proof.
\end{proof}

\bigskip

The fact that the permutation $\tau$ from Lemma~\ref{le:IItau} is an automorphism of $X$ puts some further restrictions
on parameters $m,n,r,t$ of $X \cong \Y(m,n;r,t)$. Consider the generic $8$-cycle $C$ 
from row~1 of Table~\ref{tab:classIItypes}.
Lemma~\ref{le:IItau} implies that $\tau$ maps $C$ to the $8$-cycle
from row~2 of Table~\ref{tab:classIItypes}, in particular 
$C\tau$ is of type~I. Therefore, $8$-cycles of type~I exist in $X$, and so
%%% equation
\begin{equation}\label{eq:II(r-1)^2}
(r-1)^2 = 0.
\end{equation}
%%%%
Consequently, 
$$r^m = ((r-1)+1)^m = (r-1)^m + m(r-1)^{m-1} + \cdots + m(r-1) + 1 = m(r-1) + 1,$$ and so (\ref{eq:IIr^m}) implies 
that
%%% equation
\begin{equation}\label{eq:IIrm}
 m(r-1) = 0. 
\end{equation}
%%% equation
Let $d_m$ be as in Lemma~\ref{le:IItau}. Then $\tau$
maps the directed $n$-cycle $u_0^0 u_0^1 u_0^2 \ldots u_0^{n-1}$
to $u_0^0 u_1^0 u_2^0 \ldots  u_{m-1}^{(d_m-1)t}$, and so $d_m t = 0$. Moreover, 
$d_m$ is the smallest such positive integer.
It is thus clear that 
%%% equation
\begin{equation}\label{eq:IImcongt}
|m| = |t|\ \mathrm{in}\ \ZZ_n,\ \mathrm{and}\ \mathrm{so}\ \langle m \rangle = \langle t \rangle.
\end{equation}
%%% equation
\noindent
We now show that 
%%% equation
\begin{equation}\label{eq:IImnontrivial}
|m| = |t| > 2.
\end{equation}
%%% equation
Consider the permutation $\psi$ of $V(X)$ defined by the rule: 
$u_i^j\psi = u_{i}^{-j}$ for $i \in \ZZ_m$, $j \in \ZZ_n$.
It is easy to see that $\psi$ is an automorphism of $X$
if and only if $2t = 0$. But as $\psi\rho$ interchanges adjacent vertices $u_0^0$ and $u_0^1$,
Proposition~\ref{pro:flip} implies that 
$\psi$ cannot be an automorphism of $X$, and so $2t \neq 0$, as claimed.
Note also that 
%%% equation
\begin{equation}\label{eq:IIr^2ne1}
	r^2 \notin \{\pm 1\},
\end{equation}
%%% equation
for otherwise $r^{m-i} = \pm r^i$ for all $i \in \ZZ_m$, and so 
the permutation $\varphi$ of $V(X)$
mapping according to the rule $u_0^j\varphi = u_0^{-j}$ and $u_i^j\varphi = u_{m-i}^{-j-t}$,
where $i \in \ZZ_m \setmin \{0\}$ and $j \in \ZZ_n$, is an automorphism of $X$. But this is impossible,
since then the automorphism $\varphi\rho$ inverts adjacent vertices $u_0^0$ and $u_0^1$, which contradicts
Proposition~\ref{pro:flip}.

Another immediate consequence of the existence of the automorphism $\tau$ of Lemma~\ref{le:IItau} 
is the following  lemma.

%
%   w
%

%%%%%  lemma le:IIaux  %%%%%
\begin{lemma}\label{le:IIaux}
With the notation introduced in this section let $n = md_m$, where $m$ and $d_m$ have 
the same role as in the statement of Lemma~\ref{le:IItau}. Let $a \in \{0,1,\ldots , d_m-1\}$ and
let $b \in \{0,1,\ldots , m-1\}$. Then there exists a unique $a' \in \{0,1,\ldots , d_m-1\}$
such that $(a-a')m =  (a'-a)t = b(r-1)$. 
In particular, there exists a 
unique $a' \in \{0,1,\ldots , d_m-1\}$ such that $-a'm =  a't = r-1$, and so
$m$ divides $r-1$, that is,  $r-1 \in \langle m \rangle \leq \ZZ_n$.
\end{lemma}

\begin{proof}
Lemma~\ref{le:IItau} implies that there exist unique 
$a' \in \{0,1,\ldots , d_m-1\}$ and $b' \in \{0,1,\ldots , m-1\}$ such that $am+b = (a'm + b')r$ in $\ZZ_n$. 
Let $\tau$ be as in Lemma~\ref{le:IItau}. 
Since $\tau$ maps
the outer edge $u_0^{am+b} u_1^{(a'm+b')r}$ to the inner edge
$u_b^{at} u_{b'}^{1 + a't}$, we have that $b' =  b$ and that 
$1 + a't = at + r^b$. By (\ref{eq:II(r-1)^2}) we have $r^b - 1 = b(r-1)$.
Taking into account equations (\ref{eq:IIrt=t}) and (\ref{eq:IIrm}), we see that
$(a-a')m = b(r-1)$ and $(a'-a)t = b(r-1)$, as required.
Plugging in the values $b=1$ and $a = 0$ we get that $-a'm = r-1 = a't$ in $\ZZ_n$. Since $m$ divides $n$, 
this implies that $m$ divides $r-1$, completing the proof.
\end{proof}\bigskip

We are now ready to investigate possible attachment numbers of $X$. To this end let us 
inspect the two alternating cycles containing $u_0^0$.
The directions of edges on the one on which $u_0^0$ is the tail of the two incident edges are
$$ u_0^0 \to u_1^0 \ot u_1^{-r} \to u_2^{-r} \ot 
	\cdots  \to
	u_0^{t-r-r^2-\cdots - r^{m-1}} \ot u_0^{t-1-r-\cdots -r^{m-1}} \to \cdots $$
The directions of edges on the other alternating cycle containing $u_0^0$ are
$$ u_0^0 \ot u_0^{-1} \to u_1^{-1} \ot u_1^{-1-r} \to \cdots \ot
	u_{m-1}^{-1-r-\cdots - r^{m-1}} \to u_0^{t-1-r-\cdots - r^{m-1}} \ot \cdots$$
It is therefore clear that $X$ is tightly attached if and only if there exists some 
$k \in \ZZ_n$ such that $-1 = k(t - (1+r+r^2+\cdots +r^{m-1})) - r$, that is, if and only if
\begin{equation}\label{eq:IITAcond}
r-1 = k\left(t - (1+r+r^2 + \cdots + r^{m-1})\right)\quad \mathrm{for\ some}\quad k \in \ZZ_n.
\end{equation}
Let us also note that (\ref{eq:II(r-1)^2}) implies that 
\begin{equation}\label{eq:II1+r+...}
\begin{array}{rcl} 1+r+\cdots + r^{m-1} & = & 1 + (r-1)+1 + \cdots + ((r-1)+1)^{m-1} \\
               & = & m + (1+2+\cdots + m-1)(r-1) \\
               & = & m + \frac{m(m-1)}{2}(r-1).\end{array} 
\end{equation}

%%%%%  lemma le:IIm odd  %%%%%
\begin{lemma}\label{le:IIm odd}
With the notation introduced in this section let $n = md_m$, where $m$ and $d_m$ have the same role as 
in the statement of Lemma~\ref{le:IItau}.
If either $m$ or $d_m$ is odd then $X$ is tightly attached.
\end{lemma}

\begin{proof}
We claim that in each of these two cases $ 1+r+\cdots + r^{m-1} = m$ and $\la 2(r-1)\ra = \la r-1\ra$.
Suppose first that $m$ is odd. Then (\ref{eq:IIrm}) and (\ref{eq:II1+r+...}) imply
that $ 1+r+\cdots + r^{m-1} = m$. Moreover, (\ref{eq:IIrm}) and the fact that $m$ is odd
imply that $\la 2(r-1)\ra = \la r-1\ra$.
Suppose now that $d_m$ is odd but $m$ is even. Then Lemma~\ref{le:IIaux} implies that 
$\la 2(r-1) \ra = \la r-1 \ra$. Therefore, $\frac{m}{2}(r-1) = 0$, 
and so  $ 1+r+\cdots + r^{m-1} = m$, which proves our claim.

By Lemma~\ref{le:IIaux} there 
exists a unique integer $a' \in \{0,1,\ldots , d_m - 1\}$ such that
$-a'm = a't = r-1$, and so $a'(t-m) = 2(r-1)$. 
Combining together the above claim and (\ref{eq:IITAcond}), we see that  $X$ is tightly attached.
\end{proof}

%%%%%  lemma le:II8divn  %%%%%
\begin{lemma}\label{le:II8divn}
With the notation introduced in this section let $n = md_m$, where $m$ and $d_m$ have the same role as 
in the statement of Lemma~\ref{le:IItau}.
Further, let $n = 2^i n_1$ where $n_1$ is odd. If $i \leq 2$ then $X$ is tightly attached.
\end{lemma}

\begin{proof}
By Lemma~\ref{le:IIm odd} we can assume that both $m$ and $d_m$ are even. 
Hence  $i = 2$ and $\mmod{m}{2}{4}$. 
As $m(r-1) = 0$, we either have $\frac{m}{2}(r-1) = 0$ or $\frac{m}{2}(r-1) = \frac{n}{2}$. 
We consider these two cases separately.
\smallskip

\noindent
{\sc Case~1:} $\frac{m}{2}(r-1) = 0$.\\
Then (\ref{eq:II1+r+...}) implies that $1+r+\cdots + r^{m-1} = m$. Moreover,
as $\frac{m}{2}$ is odd we have $\langle 2(r-1) \rangle = \langle r-1 \rangle$.
By Lemma~\ref{le:IIaux} there exists a unique $a' \in \{0,1,\ldots , d_m-1\}$ 
such that $-a'm = a't = r-1$, and so, as in the proof of Lemma~\ref{le:IIm odd},
$X$ is tightly attached by (\ref{eq:IITAcond}).
\smallskip

\noindent
{\sc Case~2:} $\frac{m}{2}(r-1) = \frac{n}{2}$.\\ 
Since $n$ is even, $r$ is odd, and so $r-1$ is even. Moreover,  
as $\frac{m}{2}(r-1) = \frac{n}{2}$, we have that $\mmod{r-1}{2}{4}$.
Let $c = t - m + \frac{n}{2}$, that is, $c = t - (1+r+\cdots + r^{m-1})$
in view of (\ref{eq:II1+r+...}).
By (\ref{eq:IITAcond}),  $X$ is tightly attached if and only if 
$r-1 \in \langle c \rangle$. 
By Lemma~\ref{le:IIaux}
there exists a unique $a' \in \{0,1,\ldots , d_m-1\}$ 
such that $-a'm = a't = r-1$. Note that $a'$ is odd, and so
$a'\frac{n}{2} = \frac{n}{2}$.
Consequently $a'c = 2(r-1) + \frac{n}{2}$. To show that $X$ is tightly attached it thus
suffices to see that 
$\langle 2(r-1) + \frac{n}{2} \rangle = \langle r-1 \rangle$.
It is clear that $\mmod{2(r-1) + \frac{n}{2}}{2}{4}$.
Moreover, $\frac{n}{2} \in \langle r - 1 \rangle$, and so 
$2(r-1) + \frac{n}{2} \in \langle r-1\rangle$.
Finally, let $p^t$ be any odd prime power dividing $2(r-1) + \frac{n}{2}$ and $n$.
Since $p$ is odd, $p^t$ divides $\frac{n}{2}$, and so it also divides $r-1$. Therefore,
a prime power $p^t$ dividing $n$ divides $r-1$ if and only if it divides 
$2(r-1) + \frac{n}{2}$. Since $\ZZ_n$ is a cyclic group, we indeed have 
$\langle 2(r-1) + \frac{n}{2} \rangle = \langle r-1\rangle$, as required.
\end{proof}\bigskip

Combining together the results of this section, the proof of Theorem~\ref{the:IItheorem} is now at hand.
\bigskip

%%%% PROOF OF THEOREM
\noindent
{\sc Proof of Theorem~\ref{the:IItheorem}:}\\
Throughout this proof we adopt the notation of this section.
Except for the existence of a unique 
$c \in \{0,1,\ldots , d_m - 1\}$ for which $t = cm$ and $m = ct$, claims $(i)$, $(iii)$ and $(iv)$ 
follow from (\ref{eq:IIregularly}), (\ref{eq:IIr^m}), (\ref{eq:IIrt=t}), 
Proposition~\ref{pro:IImdivn}, Lemma~\ref{le:IIaux}
and (\ref{eq:IImnontrivial}). To show that a unique such $c$ exists, observe 
that there exists a unique $c \in \{0,1,\ldots , d_m - 1\}$ for which $t = cm$. Plugging in the values 
$i = 0$ and $j = t$ in Lemma~\ref{le:IItau}
we get that $u_0^t\tau = u_0^{cm}\tau  = u_0^{ct}$. On the other hand, letting $i = 0$ and $j = m$, we get that
$u_0^m\tau = u_0^t$, and so  $u_0^t\tau = u_0^m\tau^2 = u_0^m$, which gives $m = ct$, as required.

Let us now prove $(ii)$. By contradiction, assume that 
$(\Aut X)_v$ is not isomorphic to $\ZZ_2$.
Then by Proposition~\ref{pro:IIstab},
$(\Aut X)_v \iso \ZZ_2 \times \ZZ_2$ and $m = 4$. 
Moreover, the proof of Proposition~\ref{pro:IIstab} reveals that $8$-cycles of type $V$
exist (see Table~\ref{tab:classIItypes}). In particular $3+r^2 = t$, and so 
equation (\ref{eq:II(r-1)^2}) implies that $2+2r=t$. 
Note also that (\ref{eq:II(r-1)^2}) and (\ref{eq:IIr^2ne1}) combined together 
imply  that $2(r-1) \neq 0$. By (\ref{eq:IIrm}) we thus have 
%%%%%
\begin{equation}
\label{eq:IIth2r}
	2(r-1) = \frac{n}{2}\quad \mathrm{and}\quad t = \frac{n}{2} + 4.
\end{equation}
Furthermore, Lemma~\ref{le:IIaux} and (\ref{eq:IIrm}) 
combined together imply that $n = 16n_1$ for some integer $n_1$. 
Thus either $r-1 = \frac{n}{4} = 4n_1$ or $r-1 = \frac{3n}{4} = 12n_1$. Let $a'$
be the unique element of $\{0,1,\ldots , 4n_1 -1\}$ such that $-a'm = a't = r-1$, which exists by
Lemma~\ref{le:IIaux}. As $-a'm = r-1$, we either have $a' = 3n_1$ or $a' = n_1$.
But then $a't$ equals either to $3n_1(8n_1+4) = 12n_1 + 8 n_1^2$ or to
$n_1(8n_1+4) = 4n_1 + 8n_1^2$, respectively. 
Thus, in view of the equality $r-1 = a't$, we have $8n_1^2 =   8n_1 = \frac{n}{2}$ in either case, and so
$n_1$ is odd, that is $\mmod{n}{16}{32}$. Note that this also forces $\mmod{r}{5}{8}$.

We now introduce a certain mapping $\varphi : V(X) \to V(X)$, which will be shown below to be 
an automorphism of $X$. The nature of the action of $\varphi$ will contradict
half-arc-transitivity of $X$, which thus proves that $(\Aut X)_v \cong \ZZ_2$, as claimed.
Note that (\ref{eq:II(r-1)^2}) and (\ref{eq:IIth2r}) imply that $r^2 = \frac{n}{2} + 1$, and so 
$1+r+r^2+r^3+t = 1 + r + \frac{n}{2} + 1 + \frac{n}{2} + r + t = 2t = 8$. 
Let $j \in \ZZ_n$. Then there exist unique $a \in \{0,1,\ldots , 2n_1 - 1\}$ and $b \in \{0,1,\ldots , 7\}$
such that $j = 8a + b$. The action of $\varphi$ on $u_i^j$,  $i \in \ZZ_m$, $j \in \ZZ_n$,
is given  in Table~\ref{tab:IIvarphi} and it depends on $i$ and $b$.

\begin{table}[htb]
\begin{footnotesize}
\begin{center}
\begin{tabular}{@{}c|l|l|l|l@{}}
	$b\ \backslash\ i$ & \hskip 1cm $0$ &\hskip 1cm  $1$ &\hskip 1cm $2$ &\hskip 1cm $3$ \\
	\hline  & & & & \\
   $0$ & $u_0^{-8a}$ & $u_3^{-t-8a}$ & $u_3^{-r^3-t-8a}$ & $u_2^{-r^3-t-8a}$ \\ 
   & & & &  \\
   $1$ & $u_0^{-1-8a}$ & $u_0^{-2-8a}$ & $u_3^{-2-t-8a}$ & $u_3^{-2-r^3-t-8a}$ \\ 
   & & & &  \\
   $2$ & $u_3^{-1-t-8a}$ & $u_2^{-1-t-8a}$ & $u_2^{-1-r^2-t-8a}$ & $u_1^{-1-r^2-t-8a}$ \\ 
   & & & &  \\
   $3$ & $u_3^{-1-r^3-t-8a}$ & $u_3^{-1-2r^3-t-8a}$ & $u_2^{-1-2r^3-t-8a}$ & $u_2^{-1-r^2-2r^3-t-8a}$ \\ 
   & & & &  \\
   $4$ & $u_2^{-1-r^3-t-8a}$ & $u_1^{-1-r^3-t-8a}$ & $u_1^{-1-r-r^3-t-8a}$ & $u_0^{-1-r-r^3-t-8a}$ \\ 
   & & & &  \\
   $5$ & $u_2^{-1-r^2-r^3-t-8a}$ & $u_2^{-1-2r^2-r^3-t-8a}$ & $u_1^{-1-2r^2-r^3-t-8a}$ & $u_1^{-r^2-8(a+1)}$ \\    
   & & & &  \\
   $6$ & $u_1^{-1-r^2-r^3-t-8a}$ & $u_0^{-1-r^2-r^3-t-8a}$ & $u_0^{-2-r^2-r^3-t-8a}$ & $u_3^{-2-r^2-r^3-2t-8a}$ \\
   & & & &  \\
   $7$ & $u_1^{-8(a+1)}$ & $u_1^{-r-8(a+1)}$ & $u_0^{-r-8(a+1)}$ & $u_0^{-1-r-8(a+1)}$   
\end{tabular}
\caption{The entry in $b$-th row and $i$-th column represents the image $u_i^j\varphi$ in the case
when $j = 8a + b$, where $a \in \{0,1,\ldots , 2n_1 - 1\}$ and $b \in \{0,1,\ldots , 7\}$.}
\label{tab:IIvarphi}
\end{center}
\end{footnotesize}
\end{table} 

To see that $\varphi$ is in fact a permutation of $V(X)$, we only need to observe that it is injective. 
Consider the vertices of Table~\ref{tab:IIvarphi}, which are of the form $u_0^j$. 
There are precisely eight such vertices.
It can bee seen that the congruencies modulo $8$ of their superscripts are precisely the
eight possibilities $0,1, \ldots , 7$. For instance, for $u_0^{-8a}$ we have $0$, 
for $u_0^{-1-r-r^3-t-8a}$ we have $1$, etc. 
Similarly,  one can check that there are precisely eight vertices of the 
form $u_i^j$ in Table~\ref{tab:IIvarphi} for each $i = 1,2,3$,
and that the corresponding congruencies modulo $8$
of their superscripts are again the eight possibilities $0,1, \ldots , 7$. This shows
that $\varphi$ is indeed injective and thus also bijective. It remains to be seen
that $\varphi$ preserves adjacency in  $X$. 
For the outer edges connecting $X_i$ to $X_{i+1}$, where $i \neq 3$, and for the 
inner edges of $X_0$, this is clear, as one only needs to
check that two consecutive vertices in a row or in column $0$, respectively, 
of Table~\ref{tab:IIvarphi} are adjacent.  As for the other edges, 
using the facts that $t = \frac{n}{2} + 4$, that $r^2 = \frac{n}{2} + 1$, 
and that $r \equiv 5 \pmod 8$, checking that they are indeed 
mapped to edges of $X$  is just a matter of tedious computation. We leave the details to the reader.
Therefore $\varphi$ is an automorphism of $X$.
Since it fixes $u_0^0$ and maps $u_0^1$ to $u_0^{-1}$, it follows that 
$\varphi\rho$ interchanges adjacent vertices $u_0^0$ and $u_0^{1}$ of $X$, 
which contradicts Proposition~\ref{pro:flip}.
Thus $(\Aut X)_v \cong \ZZ_2$ for all $v \in V(X)$, as claimed.

Finally, we  prove $(v)$. 
Let us suppose that $X$ is not tightly attached. 
By Lemma~\ref{le:IIm odd} $m$ and $d_m$ are both even and by Lemma~\ref{le:II8divn},
there exists some positive integer $n_1$ such that $n = 8n_1$. 
We show that $n_1 > 2$. Note first that part $(iii)$ of this theorem, 
Proposition~\ref{pro:threeorbits} and 
(\ref{eq:IImnontrivial}) combined together imply that $n_1 > 1$.
Moreover, if $n_1 = 2$, then $m = 4$, and so combining together (\ref{eq:IIrm}) and (\ref{eq:IIr^2ne1}) 
we have that $r \in \{5,13\}$.
Therefore, (\ref{eq:IImcongt}) and Lemma~\ref{le:IIaux} 
combined together imply that $t = 12 = \frac{n}{2} + 4$.
But then the mapping $\varphi$ introduced in the proof of part (ii)
is an automorphism of $X$, a contradiction. Thus $n_1 > 2$, as claimed.
Suppose now, that $n_1 > 2$ is odd and squarefree. 
Then (\ref{eq:II(r-1)^2}) implies that $r-1 = \frac{n}{2}$ (recall that $r \neq 1$). 
It follows that 
$r^2 = ((r-1)+1)^2 = (r-1)^2 + 2(r-1) + 1 = 1$, which contradicts (\ref{eq:IIr^2ne1}).
This completes the proof of Theorem~\ref{the:IItheorem}.
\hfill $\Qed$
\bigskip

There do exist connected quartic half-arc-transitive weak metacirculants of Class~II which are 
not tightly attached.
Using Theorem~\ref{the:IItheorem},  a computer search has been  performed revealing
that there are precisely $18$ such graphs of order not exceeding $1000$.
The smallest such graph is isomorphic to $\Y(4,48;13,44)$, and is of  order $192$, 
has radius $12$ and attachment number $3$. 
Table~\ref{table:IInotta} contains some information about these $18$ graphs. 
Observe that in view of existence of these graphs,  Theorem~\ref{the:IItheorem} is best possible.

\begin{table}  %%%                     TABLE OF GRAPHS OF CLASS II THAT ARE NOT TA
[htb]
\begin{center}
\begin{tabular}{||c|l|c|c||}
\hline  order & graph &  radius & att. num.\\
\hline \hline
 192  & $\Y(4,48;13,44)$  & 12 & 3 \\
 256  & $\Y(8,32;9,24)$  & 16 & 8 \\
 320  & $\Y(4,80;21,76)$  & 20 & 5 \\
 432  & $\Y(6,72;13,66)$  & 18 & 9 \\
 448  & $\Y(4,112;29,108)$ & 28 & 7 \\
 512  & $\Y(8,64;9,56)$ & 32 & 16 \\
 512  & $\Y(8,64;25,56)$  & 32 & 16 \\
 576  & $\Y(12,48;13,36)$ & 12 & 3 \\
 576  & $\Y(4,144;37,140)$ & 36 & 9 \\
 704  & $\Y(4,176;45,172)$ & 44 & 11 \\
 768  & $\Y(8,96;25,56)$ & 16 & 8 \\
 768  & $\Y(8,96;25,88)$  & 48 & 24 \\
 832  & $\Y(4,208;53,204)$ & 52 & 13 \\
 864  & $\Y(12,72;13,60)$ & 36 & 18 \\
 864  & $\Y(6,144;25,30)$  & 18 & 9 \\
 960  & $\Y(4,240;61,44)$  & 12 & 3 \\
 960  & $\Y(4,240;61,76)$  & 20 & 5 \\
 960  & $\Y(4,240;61,236)$ & 60 & 15 \\
\hline
\end{tabular}
\caption{\label{table:IInotta} All connected quartic half-arc-transitive weak metacirculants of Class~II 
of order up to $1000$ which are not tightly attached.}
\end{center}
\end{table}

Recall that a weak $(m,n)$-metacirculant is not necessarily an $(m,n)$-metacirculant. However, 
as the next proposition shows, 
quartic half-arc-transitive weak metacirculants of Class~II which are not tightly attached do have this
property.

%%%%%  proposition pro:IInotTA=>meatcirc  %%%%%
\begin{proposition}\label{pro:IInotTA=>metacirc}
Let $X$ be a connected quartic half-arc-transitive weak $(m,n)$-meta\-circulant of Class~II 
which is not tightly attached. Then $X$ is an $(m,n)$-metacirculant.
\end{proposition}

\begin{proof}
Let $r,t \in \ZZ_n$ satisfy part (iv) of Theorem~\ref{the:IItheorem}, in particular $X \cong \Y(m,n;r,t)$, 
and let the corresponding automorphisms be $\rho$ and $\sigma$. 
By part (v) of Theorem~\ref{the:IItheorem} we have that $m$ and $d_m$, where $n = md_m$, are both even
and that $n = 8n_1$,  where $n_1 > 2$ is even or not squarefree.
Observe that $X$ is a weak $(m,n)$-metacirculant relative to the ordered pair
$(\rho, \sigma\rho^k)$ for any $k \in \ZZ_n$. We show that there exists some $k \in \ZZ_n$ 
for which $\sigma\rho^k$ is of order $m$, which then completes the proof.

Let $k \in \ZZ_n$. 
Since $\sigma^m = \rho^t$, equation (\ref{eq:metagroup}) implies that
$$ (\sigma\rho^k)^m = \sigma^m\rho^{k(1+r+\cdots + r^{m-1})} = \rho^{t+k(1+r+\cdots +r^{m-1})}$$
Combining together (\ref{eq:IIrm}) and (\ref{eq:II1+r+...}),
we have two possibilities for $1+r+\cdots + r^{m-1}$. 
If $1+r+\cdots + r^{m-1} = m$,
then an appropriate $k$ exists by (\ref{eq:IImcongt}). 
We can thus assume that  $1+r+\cdots + r^{m-1} = m + \frac{n}{2}$, that is,  $\frac{m}{2}(r-1) = \frac{n}{2}$. 
Since $d_m$ is even, it is clear that $\la m + \frac{n}{2} \ra \leq \la m \ra$. 
If the order $d_0$ of $m + \frac{n}{2}$ in $\ZZ_n$ is also even, then $d_0 = d_m$. Thus
$\langle m + \frac{n}{2}\rangle = \langle m \rangle$, and we are done in view of
(\ref{eq:IImcongt}). Suppose then that $d_0$ is odd, that is $d_m = 2d_0$.
Since $m\frac{r-1}{2} = \frac{n}{2}$, we have that $\mmod{r-1}{2}{4}$. 
But then $\mmod{(r-1)^2}{4}{8}$ which contradicts (\ref{eq:II(r-1)^2}) 
and the fact that $\mmod{n}{0}{8}$.
This shows that an appropriate $k \in \ZZ_n$ does exist, as claimed.
\end{proof}\bigskip

To wrap up this section we construct an infinite family
of connected quartic half-arc-transitive weak metacirculants of Class~II 
which are not tightly attached. These graphs are constructed as regular $\ZZ_p$-covers, $p$ a prime,
of the graph $\Y(4,48;13,44)$,  the smallest example of such graphs.

%%%%% construction cons:infinite  %%%%%
\begin{construction}
\label{cons:infinite}
{\rm

Let $X$ denote the graph $\Y(4,48;13,44)$ and let $D_X$ denote the oriented graph
corresponding to the half-arc-transitive action of $\Aut X$ on $X$ in which $u_0^0 \to u_0^1$.
Let $p \geq 5$ be a prime. Following the theory developed 
in~\cite{MNS00} we construct a regular $\ZZ_p$-cover of $X$ by voltage 
assignments from the cyclic group $\ZZ_p$,  letting the
voltage of each dart corresponding to an oriented edge of $D_X$ be $1$.
Let $\C_p(X)$ denote the obtained $\ZZ_p$-cover.
Note that since $X$ is half-arc-transitive,
any automorphism of $X$ maps a cycle of $X$ with trivial net voltage to a cycle of $X$ 
with trivial net voltage. By \cite[Corollary~7.2]{MNS00}, the automorphism group $\Aut X$ lifts, that is,
there exists a group $\tilde A \leq \Aut \C_p(X)$  projecting to $\Aut X$ and acting half-arc-transitively on $\C_p(X)$. 
%Let $H = \la \rho , \sigma \ra$, where $\rho$ and $\sigma$ are as in Example~\ref{ex:Y} and let
%$\tilde{H}$ be the lift of $H$.
We now show that $\C_p(X)$ is a connected quartic half-arc-transitive weak metacirculant of Class~II with
radius $12$ and attachment number $3$. This will establish the  existence of  infinitely many 
connected quartic half-arc-transitive weak metacirculants of Class~II which are not tightly attached.

Denote the vertices of $\C_p(X)$ by 
$\{^ku_i^j\, |\, i \in \ZZ_{4},\ j\in \ZZ_{48},\ k \in \ZZ_p\}$. Let $D_{\C_p(X)}$ be the
oriented graph corresponding to the half-arc-transitive action of $\tilde A$ on 
$\C_p(X)$ such that $^0u_0^0 \to\, ^1u_0^1$.
The orientations of the edges of $D_{\C_p(X)}$ are thus
\begin{equation}
   ^ku_i^j \to \left\{\begin{array}{ccl} ^{k+1}u_{i+1}^{j} & ; & i \neq 3 \\
				    ^{k+1}u_0^{j+44} & ; & i = 3\end{array}
				    \right.\quad\quad\mathrm{and}\quad\quad
	^ku_i^j \to \ ^{k+1}u_i^{j+13^i} 				   
\end{equation}
where $i \in \ZZ_4$, $j \in \ZZ_{48}$ and $k \in \ZZ_p$.
We proceed by proving a series of claims.\smallskip

{\sc Claim~1:} There exist $r \in \ZZ_{48p}^*$ and $t \in \ZZ_{48p}$ 
such that $\C_p(X) \cong \Y(4,48p; r, t)$.\\
For any $l \in \{0,1, \ldots , 48p-1\}$ let $\alpha(l) = (j,k)$, where $j \in \{0,1,\ldots , 47\}$ and 
$k \in \{0,1,\ldots, p-1\}$ are such that
$\mmod{l}{j}{48}$ and $\mmod{l}{k}{p}$. Note that $\alpha$ is a well defined mapping. Moreover, 
since $p \geq 5$ is a prime we have $GCD(48,p) = 1$, and so $\alpha$ gives a $1-1$
correspondence of $\{0,1, \ldots , 48p-1\}$ and $\{0,1,\ldots , 47\} \times \{0,1,\ldots, p-1\}$, 
and so it also gives a $1-1$ correspondence of $\ZZ_{48p}$ and $\ZZ_{48} \times \ZZ_p$. In fact, $\alpha$
is an isomorphism of these Abelian groups.
Let $r \in \ZZ_{48p}$ be the unique element given by this
bijective correspondence such that $\alpha(r) = (13,1)$. Note that then $r \in \ZZ_{48p}^*$. 
Similarly let $t \in \ZZ_{48p}$ be the unique element such that $\alpha(t) = (44,4)$.
We now show that $\C_p(X) \cong \Y(4,48p;r, t)$. 
Let the vertex set of $\Y(4,48p;r,t)$ be $\{v_i^j\ |\ i \in \ZZ_{4}$, $j \in \ZZ_{48p}\}$, 
with edges as in Example~~\ref{ex:Y}. 
We let $\varphi : \C_p(X) \to \Y(4,48p;r,t)$ be the mapping defined by the rule 
$\varphi:\, ^ku_i^j \mapsto v_i^l$, where $l \in \ZZ_{48p}$ is the unique element such that 
$\alpha(l) = (j,k-i)$. Since $\alpha$ is
a bijection $\varphi$ is a bijection as well.
We leave the easy verification that $\varphi$ is an isomorphism of graphs
to the reader.\smallskip

{\sc Claim~2:} $\Aut \C_p(X) = \tilde A$. \\
By Claim~1 we have that $\C_p(X) \cong \Y(4,48p; r, t)$. 
It may be checked that $X$ has girth $8$. As $\C_p(X)$ is a regular $\ZZ_p$-cover of $X$
and since generic $8$-cycles of $X$ have trivial net voltage, it is clear that $\C_p(X)$ also has girth $8$. 
Moreover, any $8$-cycle of $\C_p(X)$ projects to an $8$-cycle of $X$ via the corresponding covering projection.
It thus follows that the only $8$-cycles of $\C_p(X)$ are the lifts of $8$-cycles of $X$ with net voltage $0$.
We now investigate all such $8$-cycles of $X$. 

It turns out that $X$ has precisely four $\Aut X$-orbits of $8$-cycles with voltage $0$.
These are the $\Aut X$-orbit $\HH_1$ which contains the generic $8$-cycles and the $8$-cycles of 
type~I (see Table~\ref{tab:classIItypes}),
the $\Aut X$-orbit $\HH_2$ containing the $8$-cycle 
$u_0^0u_0^1u_0^2u_3^{2-t}u_2^{2-t}u_2^{2-r^2-t}u_2^{2-2r^2-t}u_3^{2-2r^2-t}$,
the $\Aut X$-orbit $\HH_3$ containing the $8$-cycle
$u_0^0u_0^1u_1^1u_2^1u_2^{1+r^2}u_1^{1+r^2}u_1^{1-r+r^2}u_1^{1-2r+r^2}$ and
the $\Aut X$-orbit $\HH_4$ containing the $8$-cycle
$u_0^0u_0^1u_1^1u_1^{1-r}u_1^{1-2r}u_2^{1-2r}u_2^{1-2r+r^2}u_1^0$. As the covering projection
$p : \C_p(X) \to X$ gives a $1-1$ correspondence between the $8$-cycles of $\C_p(X)$ and the
$8$-cycles of $X$ with voltage $0$, we see that $\C_p(X)$ has four $\tilde A$-orbits of
$8$-cycles. They are the lifts of the $\Aut X$-orbits $\HH_i$, and so we denote the $\tilde A$-orbits of
$8$-cycles of $\C_p(X)$ with $\tilde \HH_i$ where $\HH_i = p(\tilde \HH_i)$.

Observe that $\C_p(X)$ has precisely four $\tilde A$-orbits of $2$-paths (where no
distinction is made on the orientation of these paths).
Let $\PP_1$ be the orbit containing the $2$-path $^0u_0^0\, ^1u_0^1\, ^2u_0^2$, let $\PP_2$
be the orbit containing the $2$-path $^0u_0^0\, ^1u_0^1\, ^2u_1^1$,
let $\PP_3$ be the orbit containing the $2$-path $^0u_0^0\, ^1u_0^1\, ^0u_3^5$ and let
$\PP_4$ be the orbit containing the $2$-path $^1u_0^1\, ^0u_0^0\, ^1u_1^0$.
We now consider the bipartite graph $Bip_{2,8}$ whose vertex set
is the union of the set of $8$-cycles of $\C_p(X)$ and the set of $2$-paths of $\C_p(X)$ with
a $2$-path $P$ being adjacent to an $8$-cycle $C$ if and only if $C$ contains $P$.
It is straightforward to check that the valency in $Bip_{2,8}$ of any $2$-path 
from any one of $\PP_2$, $\PP_3$ and $\PP_4$ is $8$, whereas the valency
of any $2$-path from $\PP_1$ is $4$.

Suppose now that $\C_p(X)$ is arc-transitive. 
Following \cite{DM05}, we assign  a letter $D$, $A^+$ or $A^-$ to an internal vertex $v$ 
of a $2$-path $P$  in $\C_p(X)$ depending on whether 
$v$ is the head of one and tail of the other, head of both, or tail of both 
of the two incident edges in $P$, respectively. 
In such a way a code of length $2$ can be assigned to each $3$-path in $\C_p(X)$. 
For example, a directed $3$-path is assigned the code $D^2$.
Now, by \cite[Lemma~2.1]{DM05}
either every automorphism of $\C_p(X)$ preserves or reverses the orientation of every edge or
any two codes of length $2$ are permutable  (by an automorphism of $\C_p(X)$).
If the latter occurs then there exists an automorphism of $\C_p(X)$ 
mapping some $3$-path with code $DA^+$ to
some $3$-path with code $A^-A^+$. 
It can be seen that any $3$-path with code $DA^+$ lies on some $8$-cycle from $\tilde \HH_4$,
whereas the only $8$-cycles of $\C_p(X)$ which contain $3$-paths with code $A^-A^+$ are the
$8$-cycles from $\tilde \HH_1$. Therefore, some automorphism of $\C_p(X)$, mapping an $8$-cycle 
from $\tilde \HH_4$ to an $8$-cycle from $\tilde \HH_1$ exists. However, it may be verified that
$8$-cycles from $\tilde \HH_4$ contain two $2$-paths from $\PP_1$, whereas $8$-cycles from
$\tilde \HH_1$ contain none, which contradicts the fact that the $2$-paths from $\PP_1$ are the only
$2$-paths of $\C_p(X)$ of valency $4$ in $Bip_{2,8}$.
Therefore, either every automorphism of $\C_p(X)$  preserves  the orientation of every edge of $\C_p(X)$
or every automorphism of $\C_p(X)$  reverses  the orientation of every edge of $\C_p(X)$.
More precisely, by assumption about arc-transitivity of $\C_p(X)$,
there must exist an automorphism reversing the  edges orientations.
Let $\tilde G$ be the subgroup of index $2$ of  $\Aut \C_p(X)$ 
consisting of all the automorphisms preserving the orientation of the edges. 
If $\tilde A = \tilde G$, then  $\tilde A$ is an index $2$ subgroup of $\Aut \C_p(X)$ and is thus normal.
Since $|\Aut X| = 2\cdot 4 \cdot 48 = 2^{7}\cdot 3$, we have that $|\tilde A| = p\cdot 2^7\cdot 3$, and so the group of covering
transformations $CT_X = \ZZ_p$ is the unique p-Sylow subgroup of $\tilde A$. It is now clear that $CT_X$ is normal in $\Aut \C_p(X)$,
and so $\Aut \C_p(X)$ projects. Thus $X$ is arc-transitive, a contradiction. 
Therefore $|\tilde G_{^0u_0^0}| \geq 4$. However, this is also impossible as then
$\PP_1$ and $\PP_2$ are contained in the same $\tilde G$-orbit which is clearly impossible in
view of the valencies of the respective $2$-paths in $Bip_{2,8}$.
This contradiction finally proves that $\Aut \C_p(X) = \tilde A$, and so $\C_p(X)$ is 
half-arc-transitive, as claimed.\medskip

{\sc Claim~3:} $Y_p$ has radius $12$ and attachment number $3$. \\
By Claim~2 the orientation of the edges of $\C_p(X)$ given by $\tilde A$ 
is in fact an orientation given by the half-arc-transitive action of $\Aut \C_p(X)$, and so 
it is clear that an alternating cycle of
$X$ is lifted into an alternating cycle of $\C_p(X)$. Therefore, the radius of $\C_p(X)$ is $12$.
The fact that the attachment number is $3$ is now also clear.
\medskip

Combining together Claims~1, 2 and 3 we thus see that
$\C_p(X)$ is a connected quartic half-arc-transitive weak metacirculant of Class~II
with radius $12$ and attachment number $3$, as claimed. Note that the graph $\Y(4,240;61,44)$ from
Table~\ref{table:IInotta} is isomorphic to $\C_5(X)$.
}
\end{construction}

%%%%%%%%%%%%%%%%%%%%%%%%%%%%%%%%%%%%%%%%%%%%%%%%%%%%
%%%%%%%%%%%%%%%%%%%%%%%%%%%%%%%%%%%%%%%%%%%%%%%%%%%%
%%%%%%%%%%%%%%  Conclusions  %%%%%%%%%%%%%%%%%%%%%%%
%%%%%%%%%%%%%%%%%%%%%%%%%%%%%%%%%%%%%%%%%%%%%%%%%%%%
%%%%%%%%%%%%%%%%%%%%%%%%%%%%%%%%%%%%%%%%%%%%%%%%%%%%

\section{Conclusions}
\label{sec:conclusions}

\indent

Half-arc-transitive metacirculants of Classes III and IV will be studied 
in a sequel to this paper, together with their connection to graphs in Classes I and II.

Let us mention, however, the existence of an infinite 
family of connected quartic half-arc-transitive metacirculants of Class~IV, which
are loosely attached, and thus not in Class~I. 
They arise as certain $\ZZ_p$-covers, $p \geq 7$  a prime, in the following way.
We start with $X = \Z(20,5;9,2)$, a loosely attached half-arc-transitive metacirculant of Class~IV,
as our base graph. (Recall that the graphs $\Z(m,n;r,t)$ were defined in Section~\ref{sec:examples}.)
Of the two possible $\Aut X$-admissible orientations of the edges of $X$,
we choose the one in which  $u_0^0$ is the tail of the edges $u_0^0u_1^0$ and $u_0^0u_9^0$, and 
denote the corresponding oriented graph by $D_X$. 
Next, to every arc of $X$ assign voltage $1$ or $-1$ in $\ZZ_p$ depending on whether its orientation
is or is not compatible with the orientation of the corresponding edge in $D_X$.
Call the obtained covering graph $\C_p(X)$.
Using similar techniques as in Construction~\ref{cons:infinite}
one can show that  $\C_p(X)$  is isomorphic to $\Z(20p,5;k,2)$, where
$k \in \ZZ_{20p}$ is the unique element such that $\mmod{k}{9}{20}$ and $\mmod{k}{1}{p}$. Moreover,
one can also see that $\C_p(X)$ is a loosely attached half-arc-transitive graph, thus giving
an infinite family of half-arc-transitive metacirculants of Class~IV which are not tightly attached.
The technical details are omitted.

%%%%%%%%%%%%%%%%%%%%%%%%%%%%%%%%%%%%%%%%%%%%%%%%%%%%
%%%%%%%%%%%%%%%%%%%%%%%%%%%%%%%%%%%%%%%%%%%%%%%%%%%%
%%%%%%%%%%%%%  Bibliography  %%%%%%%%%%%%%%%%%%%%%%%
%%%%%%%%%%%%%%%%%%%%%%%%%%%%%%%%%%%%%%%%%%%%%%%%%%%%
%%%%%%%%%%%%%%%%%%%%%%%%%%%%%%%%%%%%%%%%%%%%%%%%%%%%

\end{document}